\documentclass[11pt,letterpaper]{amsart}
\usepackage{amssymb, amsthm, mathtools,bbm}
\usepackage{amsaddr}
\usepackage[square,sort,comma,numbers]{natbib}
\usepackage[none]{hyphenat}
\newtheorem{theorem}{Theorem}
\newtheorem{corollary}[theorem]{Corollary}
\newtheorem{lemma}[theorem]{Lemma}
\newtheorem{conjecture}[theorem]{Conjecture}
\newtheorem{problem}[theorem]{Open problem}
\newtheorem{definition}[theorem]{Definition}
\newtheorem{proposition}[theorem]{Proposition}

\numberwithin{equation}{section}
\numberwithin{theorem}{section}

\newcommand{\rr}{{\mathbb{R}}}

\newcommand{\zz}{{\mathbb{Z}}}
\newcommand{\nn}{{\mathbb{N}}}

\newcommand{\ee}{{\mathbb{E}\,}}

\newcommand{\pp}{{\mathbb{P}}}

\newcommand{\beq}[1]{\begin{equation} \label{#1}}
\newcommand{\eeq}{\end{equation}}

\newcommand{\lil}{\mathrm{LIL}}
\newcommand{\normal}{\mathcal{N}}
\newcommand{\F}{\mathcal{F}}

\begin{document}

	\title[On Transcendentally Normal Numbers]{ Transcendence Meets Normality: Construction of Transcendentally Normal Numbers}
     \author{Chokri Manai}
     \address{Courant Institute of Mathematical Sciences, New York University}
     \email{cm7411@nyu.edu}
        \date{\today}
	
	\begin{abstract}
	In this work, we study real numbers $x$ for which $p(x)$ is (absolutely) normal for every non-constant polynomial with integer coefficients $p$. We call such numbers  \textit{transcendentally normal}. We prove that almost every real number is transcendentally normal and provide an explicit construction of such a number, based on Sierpiński’s covering method and  novel ideas involving stretch functions and P\'{o}lya-type inequalities. In the next step, we transform this construction into an algorithm that computes the digits of a t-normal number recursively in all integer bases. Moreover, we extend our covering approach to construct and compute LIL-normal numbers whose discrepancies are of the order predicted by the law of the iterated logarithm. We also take the opportunity to discuss several interesting open problems.
	\end{abstract}
	\maketitle

\section{Introduction}\label{sec:intro}

Since Borel showed that almost all real numbers are (absolutely) normal \cite{Bor09}, one of the biggest challenges in the field has been to prove that specific mathematical constants such as $e$ and $\pi$ are normal. This problem remains wide open. In fact, constructing even a single normal number is surprisingly difficult; this was first achieved independently by Lebesgue and Sierpiński \cite{Leb17, Sir17}. The first algorithm was proposed by Alan Turing \cite{Tur92},  and his ideas have been recently made rigorous in  \cite{BFP07}. Today, there are several efficient algorithms for the computation of normal numbers \cite{ABSS17, BF02, BHS13, BY19, AB11, LM21}; the fastest currently known construction runs in nearly linear time \cite{LM21}. In this paper, we want to extend the construction problem by a significant degree of difficulty. To motivate our problem we propose the reasonable conjecture that not only Euler's number $e$ itself is a normal number, but also $p(e)$ for any non-constant polynomial with integer coefficients $p$. We call such numbers \textit{transcendentally normal}, or simply  t-normal. In this context, the natural question arises: is there an explicit example of a transcendentally normal number? Our first main result in this paper is an affirmative answer and the specification of an algorithm computing a t-normal number.

\begin{theorem}\label{thm:main}
  The set $\Omega_{\infty}$  of transcendentally normal numbers is of full measure. Moreover, the explicit "Sierpiński numbers" $(\xi_{\mathcal{P}}(r))_{r \in (0,1]}$ \textnormal{(}see \eqref{eq:sirnumP}\textnormal{)} are transcendentally normal. The algorithm proposed in Section~\ref{sec:alg} computes a transcendentally normal number $\nu \in (0,1) $ digit-by-digit in all integer bases $b \geq 2$. In particular $\nu$ is computable.
\end{theorem}
\noindent In fact, in Section~\ref{sec:const} we consider more general countable families of continuous functions $\mathcal{F}$ and under mild conditions on $\mathcal{F}$, we give a construction of Sierpiński numbers $(\xi_{\mathcal{F}}(r))_{r \in (0,1]}$  for which all numbers $(f(\xi_{\mathcal{F}}(r)))_{f \in \mathcal{F}}$ are simultaneously normal. The corresponding Theorem~\ref{thm:cover} can be seen as an (almost) optimal result on simultaneous normality.  The algorithm computing $\nu$ in Section~\ref{sec:alg} relies on some specific quantitative bounds on polynomials (see Proposition~\ref{prop:poly}) and thus does not address more general countable families of functions $\mathcal{F}$.  However, for a given family of (regular enough) functions $\mathcal{F}$  the algorithmic approach can be easily modified.  As far as we know, a similar problem has only been discussed in \cite{BM22}, where the authors present a continued-fraction type algorithm for the computation of a real number $x$ such that $x$ and its reciprocal $1/x$ are normal. Not only do we propose a different approach based on Sierpiński's cover method \cite{Sir17, BF02}, but we also vastly generalize the results in \cite{BM22}. We stress that prior to this work, even a concrete example for a number $x$ such that $x$ and $x^2$ are normal has not been known.

It is often of interest how quickly the digit frequencies of a normal number converge to their asymptotic value. This can be addressed in terms of the discrepancy number $D_N$. For a sequence of real numbers $(x_j)_{j \in \nn} \subset [0,1) $  one defines \cite{ABSS17}
\begin{equation}
    D_N((x_j)_j) := \sup_{0 \leq \alpha < \beta \leq 1} \left| \frac1N \sum_{j =1}^{N} \mathbbm{1}_{[\alpha,\beta)}(x_j) - (\beta - \alpha) \right|, 
\end{equation}
which measures how equidistributed the numbers $x_1, \ldots, x_N $ are. For any base $b \geq 2$ (or more generally any real number $b > 1$) we introduce the function
\begin{equation}\label{eq:disc}
D_N^{b}(x) := D_N((\{b^{j-1}x\})_j),
\end{equation}
where $\{ \cdot \}$ denotes the fractional part of a real number. A real number $x$ is normal if and only if $\lim_{N \to \infty} D_N^{b}(x) = 0$ for all bases $b \geq 2$. Note that the discrepancy number $D_N^{b}(x)$ slickly takes into account the frequency of all finite $b$-ary words in a single function, and thus it is desirable to bound $D_N^{b}(x)$ for some given (normal) number. Unfortunately, the rather elementary algorithm computing the number $\nu$ in Theorem~\ref{thm:main} does not yield any control on the discrepancy number.  Our second result will remove this weakness. To motivate our findings, we first recall the behavior of $D_N^b(x)$ for a random number $x \in [0,1]$. In a series of works \cite{Fu08,Phi75,GG64}  the following version of the law of the iterated logarithm (LIL) has been established as almost sure identity,
\begin{equation}\label{eq:LIL}
    \limsup_{N \to \infty} \frac{\sqrt{N} D_N^{b}(x)}{\sqrt{\log \log N }} = L_b := \begin{cases}
        \sqrt{\frac{84}{81}},& \text{ if } \, b = 2, \\
        \sqrt{\frac{b+1}{2(b-1)}}, & \text{ if } \, b \text{ odd} , \\
        \sqrt{\frac{b(b+1)(b-2)}{2(b-1)^3}} & \text{ if } \, b \geq 4 \text{ even.} 
    \end{cases}
\end{equation}
Although \eqref{eq:LIL} is reminiscent of the LIL in the context of simple random walks and Brownian motions, its proof is more involved and requires a distinct analysis. LIL-type estimates for the discrepancy number were originally conjectured by Erd\H{o}s and G\'{a}l and the full solution \eqref{eq:LIL} for all bases $b \geq 2$ (in fact, for all positive real numbers $\theta >1$) was given by Fukuyama 40 years later \cite{Fu08}.  We call a number $x$ for which \eqref{eq:LIL} holds for all bases $b \geq 2$ an \textit{LIL-normal} number and analogously one may also define \textit{transcendentally LIL-normal numbers}. Our second result establishes their computability.
\begin{theorem}\label{thm:main2}
  The set $\Omega_{\infty}^{\lil}$  of transcendentally LIL-normal numbers is of full measure. Moreover, the modified Sierpiński numbers $(\xi^{\lil}_{\mathcal{P}}(r))_{r \in (0,1]}$ \textnormal{(}see \eqref{eq:sirnumPLIL}\textnormal{)} are transcendentally LIL-normal. The algorithm proposed in Section~\ref{sec:alg} computes a transcendentally LIL-normal number $\nu_{\lil} \in (0,1) $ digit by digit in all bases $b \geq 2$. In particular $\nu_\lil$ is computable.  
\end{theorem}

Some comments are in order.
\begin{enumerate}
    \item Again, one may in fact allow for more general countable families $\mathcal{F}$ for which the modified Sierpiński numbers are simultaneously LIL-normal (see Corollary~\ref{cor:FLIL}).
    \item Note that the LIL for the discrepancy number $D_N^{b}(x)$ implies LIL-type bound for the frequencies of all finite $b$-ary strings. For example, \eqref{eq:LIL} guarantees that the number of occurrences of a fixed digit among the first $N$ base-$b$ digits of $\nu_\lil$ deviates at most $\mathcal{O}(\sqrt{N \log \log N})$ from the expected value $N/b$.
    \item In some sense, $\nu_\lil$ can be considered as a "fully generic number" - even though strictly speaking such a number does not exist \cite{Knu97}. The digits of $\nu_\lil$ behave exactly as one would expect for a randomly chosen number and that remains even true after applying elementary algebraic operations. In that regard, Theorem~\ref{thm:main2} is optimal.
    \item Our proof yields non-asymptotic bounds on the discrepancy number of the form
    $$ D_N^b(p(\nu)) \leq C(b,p) \sqrt{\frac{\log \log N}{N}} $$
    for any $p \in \zz[X]$ with some computable constant $C(b,p)$. Even if one only seeks for the best possible algorithm in terms of the discrepancy number $D_N^b$ and that only for the case $ p = x$ - at the possible expense of sacrificing some of the number’s random character - the above bound is essentially state of the art and beats for example Levin's remarkable approach \cite{AB11, Lev79}. We only know about the rather recent construction from \cite{ABSS17}, which is slightly better in this regard, as it produces a bound of the form $D_N^b(\nu) \leq C(b)N^{-1/2}$. It is a very difficult, yet also interesting, problem to construct numbers with lower discrepancy numbers. One should note that even the existence of numbers with lower discrepancy numbers is not guaranteed and an open problem due to Bugeaud and Korobov asks for the best possible bound on the discrepancy numbers.
\end{enumerate}

So far, we have not discussed the complexity of our algorithms for the simple reason that the algorithms in this work have exponential or worse running time and are thus by no means efficient. On the other hand, it is very unclear to us whether efficient algorithms should be expected for the computational tasks considered. We propose that as a potentially stimulating problem.
\begin{problem}
    Is there an efficient / polynomial algorithm computing an LIL-normal number / a transcendentally normal number?
\end{problem}

Besides the computational interest, we feel that the combination of algebraic concepts and normality is of inherent mathematical importance. This was in fact the original motivation for this work and we would like to discuss that background in the following. Let us first recall two types of results that have been established so far for normal numbers. First, it has been shown that that the geometry of the set of non-normal numbers is very intricate and in particular non-normal numbers form a comeager set \cite{AL23, Man25, Ol04a, Ol04b}. Secondly, there is a range of beautiful results regarding normality properties for different bases \cite{Bug12, Cas59, Egg49, Schm62}. Nevertheless, the behavior of normality under algebraic operations is largely not understood.  Two notable exceptions are Rauzy's theorem on the so-called deterministic numbers - which preserve normality when added to a normal number \cite{Rau76} - and the recent result by Dayan--Ganguly--Weiss \cite{DGW24} who show that only non-zero rational multipliers preserve normality. An overview of normality preserving operations is given in \cite{Man26b}. Unfortunately, a proof of the following conjecture seems elusive at the moment.  

\begin{conjecture}\label{con:alg}
    Every irrational algebraic number $\alpha \in \rr$ is normal.
\end{conjecture}

We are currently very far from proving Conjecture~\ref{con:alg}: so far, it has only been shown that the number of '1's within the first $N$ digits in the binary representation of irrational algebraic numbers grows at least as fast as $\sqrt{N}$ \cite{BBS16,Van18}. Most likely irrational algebraic numbers are normal for very different reasons than transcendent constants such as $e$ and $\pi$. The concept of transcendental normality sets $e$ and $\pi$ nicely apart from algebraic numbers. More generally, one may introduce an analogue of the algebraic degree in the context of normal numbers. 

\begin{definition}\label{def:deg}
    For a real number $x \in \rr$ we define its algebraic normality degree as
    \begin{align}
        \deg_{an}(x) := \inf \{ k \in \nn  \, | \, &\text{ there exists a  polynomial } p \in \zz[X] \\ &\text{ with } \deg(p) = k \text{ s.t. } p(x) \text{ is not normal} \}
    \end{align}
    if the set on the RHS is not empty, and otherwise $\deg_{an}(x) = \infty$. A number $x$ is transcendentally normal if and only if $\deg_{an}(x) = \infty.$  
\end{definition}

If we denote by $\deg_a$ the ordinary algebraic degree of a number, we obviously have $\deg_{an}(x) \leq \deg_a(x)$. In particular a transcendentally normal $x$ is transcendent in the usual sense. In the above definition,  considering rational polynomials $p\in \mathbb{Q}[X]$ instead of integer polynomials does not affect the algebraic normality degree. If one further allows the coefficients of $p$ to be algebraic numbers $\alpha \in \overline{\mathbb{Q}}$, a finite algebraic normality degree might be lowered in view of Conjecture~\ref{con:alg}. However, an interesting open problem is to prove or disprove whether $p(x)$ is normal for any t-normal number $x$ and any algebraic non-constant polynomial $p \in \overline{\mathbb{Q}}[X]$ with algebraic coefficients.

The concept of normality degree allows a reformulation of Conjecture~\ref{con:alg}. Namely, Conjecture~\ref{con:alg} is equivalent to the
equality 
\begin{equation}
    \deg_{an}(x) = \deg_a(x)  \quad \text{ if } \deg_a(x) < \infty.
\end{equation}
This motivates the study of the geometry of the sets $\Omega_k$ of numbers with algebraic normality degree $k \in \nn \cup \{\infty \}$ to gain some insight on Conjecture~\ref{con:alg}. At this point, the methods in the field are not strong enough to show anything of interest about the sets $\Omega_k$. It is not even known whether $\Omega_k \neq \emptyset$ for $2 \leq k \in \nn. $ So far, only some recent partial results exists which imply the existence of a non-normal number $x$ for which $\varphi(x)$ is normal provided $\varphi$ is a non-affine polynomial or, more generally, an analytic function \cite{ABS22, BB25, HS15, Man26}.  This illustrates how currently unattainable a proof of Conjecture~\ref{con:alg} appears to be.  We propose as a first step towards a better understanding of the behavior of normal numbers under algebraic operations the following interesting question:
\begin{problem}
    Show that the sets $\Omega_k$ are of Hausdorff dimension 1 (or at least of positive Hausdorff dimension). Construct an explicit number $x \in \Omega_k$ for all (some) $k \in \nn$.
\end{problem}

In this work, we begin this mathematical journey by studying the set of t-normal numbers $\Omega_{\infty}$, which is of full measure in contrast to the sets of finite algebraic normality degree (see Corollary~\ref{cor:full} below).  Hence, the existence of a t-normal number is guaranteed and the real challenge here is to give an explicit example of a t-normal number, i.e. to "find a hay in the haystack". This is exactly the purpose of our main results Theorem~\ref{thm:main} and Theorem~\ref{thm:main2}. In the course of proving our results, we provide techniques that are robust enough to construct a simultaneously normal number for more general countable families of functions.

We briefly sketch the structure of this work and give an overview of our methods:

\begin{enumerate}
    \item In Section~\ref{sec:pre}, we introduce some notation  and show that for any countable collection of nonsingular functions $(f)_{f \in \F}$  all numbers $f(x)$ are simultaneously normal for almost every $x \in \rr$ (see Prop.~\ref{prop:full}). A key technical tool, the so-called stretch function $s_f$ associated to a measurable function $f$, is defined (see Definition~\ref{def:stretch}). Most importantly we show that the stretch function $s_f$ is continuous (see Lemma~\ref{lem:stretch}) and derive a convenient bound in the case of polynomials with integer coefficients, which depends only on the degree of the polynomial (see Proposition~\ref{prop:poly}).
    \item In Section~\ref{sec:const}, we generalize Sierpiński's argument \cite{Sir17} to give an explicit construction of $\F$-normal (and more specifically of a t-normal) numbers. We first recall Sierpiński's construction in the original setup of normal numbers following the modernized presentation \cite{BF02}. In a next step, we explain how the tools of Section~\ref{sec:pre} can be used to construct modified covers $U_\F(r)$ coping with countable families of functions $\F$ and corresponding generalized Sierpiński numbers $\xi_\F(r)$, which turn out to be $\F$-normal (see Theorem~\ref{thm:cover}). We close this section by discussing some interesting properties of the Sierpiński numbers $\xi_\F(r)$ (see Proposition~\ref{prop:sirprop} and Proposition~\ref{prop:comenu}).
    \item In Section~\ref{sec:lil} we extend our Sierpiński-type approach from Section~\ref{sec:const} to deal with the case of LIL-normality.  Section~\ref{sec:lil} is the most technical part of this work. The main challenge is not the simultaneous control for all non-constant polynomials with integer coefficients - this can be addressed exactly as in Section~\ref{sec:const} before -but rather the replacement of Sierpiński's cover by a finer analysis up to the LIL case. This is largely due to the fact that the proof of the law of iterated logarithm for the discrepancy number is technically involved and constructing a LIL-cover essentially requires turning all the asymptotic assertions in the proof to non-asymptotic bounds. Since Fukuyama's proof relies on many previous works \cite{Phi75, Str67, Ber76} it is rather tedious to use exactly his argument as a blueprint for our construction. Instead, we give a simplified approach which is still inspired by the works \cite{Ber76, Fu08} but allows for a self-contained approach without martingale comparisons and invariance principles.
    \item In the last Section~\ref{sec:alg}, we deal with the computability of transcendentally (LIL)-normal numbers.  Based on ideas in \cite{BF02}, we turn the constructions of Section~\ref{sec:const} and Section~\ref{sec:lil} into an algorithm. We first discuss the theoretical backbone of our algorithm in form of Proposition~\ref{prop:algo} which works under the assumption that our covers are weakly computably measurable (see Definition~\ref{def:meas}). We believe that this approach is in fact flexible enough to deal with  many "hay in the haystack" problems, where elements with a desired property are highly likely with respect to some probability measure, but it is hard to give examples for such elements.  In a second step, we show how to cope with the countable collection of polynomials: the preimages of polynomials preserve the computable measurability under the assumption that truncations of countable unions can be controlled in a computable manner (see Proposition~\ref{prop:polywcm}). This part largely relies on Proposition~\ref{prop:poly}. In a final step, we put the pieces together and give explicit estimates for the truncations of our Sierpiński-type covers. 
\end{enumerate}

\section{Preliminary Considerations}\label{sec:pre}
In the following, $\mu$ and $\mu^*$ denote Lebesgue measure and Lebesgue outer measure on the real line, respectively. Integrals and null sets are understood with respect to $\mu$, and measurability always means Lebesgue measurability. Moreover, we recall some standard terminology. We say that a real number $x \in \rr$ is simply $b$-normal for some integer $b \geq 2$ if all digits $\{0,1,\ldots, b-1\}$ are asymptotically equidistributed in the $b$-ary representation of $x$. If not just the digits but all finite words $w \in \bigcup_{k = 1}^{\infty}\{0,1,\ldots, b-1\}^k$ appear with asymptotic frequency $b^{-|w|}$ in the $b$-ary decomposition, we call $x$ a $b$-normal number. Finally, $x$ is an (absolutely) normal number if $x$ is $b$-normal with respect to every integer base $b \geq 2$. We denote by $\mathcal{N}$ the set of normal numbers, and $\mathcal{N}^c$ stands for its complement, the set of non-normal numbers.
We recall the nonstandard notion of LIL-normality for real numbers satisfying the law of the iterated logarithm \eqref{eq:LIL}. The respective set is denoted by $\mathcal{N}^{\lil}$ and of course $\mathcal{N}^{\lil} \subset \mathcal{N}.$

It is natural to generalize the notion of t-normal numbers to an arbitrary collection of functions.
\begin{definition}
    Let $\mathcal{F} \subset \rr^{\rr}$ be a collection of real-valued functions. We say that a number $x \in \rr$ is $\mathcal{F}$-normal if $f(x)$ is a normal number for every $f \in \mathcal{F}$. The set of $\mathcal{F}$-normal numbers is denoted by $\mathcal{N}_{\mathcal{F}}.$ Similarly, $x \in \rr$ is $\mathcal{F}^\lil$-normal if $f(x)$ is a LIL-normal number for every $f \in \mathcal{F}$. The $\mathcal{F}^\lil$-normal numbers form the set $\normal_{\F}^{\lil}$. 
\end{definition}

 In the next two subsections, we characterize the countable sets $\mathcal{F}$ for which almost every number is  $\mathcal{F}$-normal (or $\mathcal{F}^\lil$-normal) and introduce some technical tools on which our construction of t-normal numbers will rely.
 
\subsection{Nonsingular functions}\label{sec:nonsing}

  We recall the concept of non-singularity for real-valued functions.

\begin{definition}\label{def:ns}
    A function $\varphi : \rr \to \rr$ is called nonsingular if, for every null set $N \subset \rr$, one has $\mu^*(\varphi^{-1}(N))=0$.
\end{definition}

Note that we do not require nonsingular functions to be measurable. A set of outer measure zero is Lebesgue measurable, so the definition is compatible with all later measure statements. Non-singularity should not be confused with Lusin's $N$-property, which requires $\mu^*(\varphi(N))=0$ for every null set $N \subset \rr$. The following sufficient criterion for non-singularity is folklore.

\begin{lemma}~\label{lem:smooth}
    Let $\varphi : \rr \to \rr$ be a $C^1$-function. Then, if 
    \begin{equation}
        \mu(\{ x \in \rr \, | \, \varphi'(x) = 0 \}) = 0,
    \end{equation}
    the function $\varphi$ is nonsingular. In particular, every non-constant polynomial $p$ is nonsingular.
\end{lemma}

To be self-contained, we sketch the easy proof of Lemma~\ref{lem:smooth}.
\begin{proof}
Let $N \subset \rr$ be a null set. For $m \in \nn$ consider the open set
\[
 U_m:=\{x\in\rr:|\varphi'(x)|>1/m\}=\bigcup_{k=1}^{\infty}I_{k,m},
\]
where the $I_{k,m}$ are its pairwise disjoint connected components. Since $\varphi'$ has a constant sign on each $I_{k,m}$, the restriction of $\varphi$ to $I_{k,m}$ is injective. The change-of-variables formula therefore gives
\begin{align*}
 \mu\bigl(\varphi^{-1}(N)\cap I_{k,m}\bigr)
 &=\int_{\varphi(I_{k,m})}\mathbbm{1}_{N}(y)
   \frac{dy}{|\varphi'(\varphi^{-1}(y))|}
 \leq m\mu(N)=0.
\end{align*}
Moreover,
\[
 \{x:\varphi'(x)\neq0\}=\bigcup_{m=1}^{\infty}U_m.
\]
Consequently,
\begin{align*}
 \mu(\varphi^{-1}(N))
 &\leq \mu(\{x:\varphi'(x)=0\})
   +\sum_{m=1}^{\infty}\sum_{k=1}^{\infty}
     \mu\bigl(\varphi^{-1}(N)\cap I_{k,m}\bigr)=0.
\end{align*}
\end{proof}
Note that Lemma~\ref{lem:smooth} implies that every non-constant polynomial $p$ is nonsingular since the set $\{ x \in \rr \, | \, p'(x) = 0 \}$ is always finite. 

The following proposition gives a simple sufficient criterion for the set of $\mathcal{F}^\lil$-normal numbers to be of full measure.
\begin{proposition}\label{prop:full}
    Let $\mathcal{F}$ be a countable set of real-valued functions. If every $f \in \mathcal{F}$ is nonsingular, 
    the set of $\mathcal{F}^\lil$-normal numbers is of full measure, i.e.
        \begin{equation}
            \mu(\{ x \in \rr \, | \, x \text{ is not } \mathcal{F}^\lil\text{-normal} \} ) = 0.
        \end{equation}
   
\end{proposition}

\begin{proof}
    We observe that the set 
    \begin{equation}
        (\mathcal{N}_{\mathcal{F}}^\lil)^c = \bigcup_{f \in \mathcal{F}} f^{-1}((\mathcal{N}^\lil)^c)
    \end{equation}
    is the countable union of null sets, since  $(\mathcal{N}^\lil)^c$ has zero measure in view of the law of the iterated logarithm \cite{Fu08}  and each $f \in \mathcal{F}$ is nonsingular. Thus, $\mu((\mathcal{N}_{\mathcal{F}}^\lil)^c) = 0.$
\end{proof}
Note that if one  only wants to show that the $\F$-normal numbers form a set of full measure, one only needs to employ Borel's theorem \cite{Bor09} and the proof of Proposition~\ref{prop:full} becomes completely elementary. Of course, Proposition~\ref{prop:full} is not a deep result, yet it is rather general and useful for our purposes. In particular, we arrive at the following
\begin{corollary}\label{cor:full}
     The set $\Omega_{\infty}$ of transcendentally normal numbers and the set $\Omega_{\infty}^\lil$ of transcendentally LIL-normal numbers are of full measure, but meager (or of first category).
\end{corollary}
\begin{proof}
    Since non-constant polynomials are nonsingular, $\Omega_{\infty}^\lil$ (and $\Omega_{\infty}$ is a larger set) is of full measure by Proposition~\ref{prop:full}. Since the normal numbers form a meager set \cite{Ol04a,Ol04b} and taking $p(X)=X$ shows that every t-normal number is normal, both $\Omega_{\infty}$ and $\Omega_{\infty}^{\lil}$ are subsets of a meager set.
\end{proof}

\subsection{Stretch functions: general properties}
While non-singularity gives an easy condition for $\mathcal{F}$-normal numbers to be generic, its qualitative nature does not allow for an explicit construction. To obtain quantitative control on how a map distorts the size of sets, we introduce the stretch function $s_f$ for measurable functions $f$.

\begin{definition}\label{def:stretch}
Let $I=[0,1]$ and let $f:I\to\rr$ be measurable. We define the corresponding stretch function by
\begin{equation}\label{eq:stretch}
 s_f(\delta):=\sup\{\mu(f^{-1}(A)):A\subset\rr\text{ measurable and }\mu(A)\leq\delta\}.
\end{equation}
\end{definition}

Let us collect some elementary properties of the stretch function.
\begin{lemma}\label{lem:stretch}
Let $f:[0,1]\to\rr$ be measurable with stretch function $s_f$.
\begin{enumerate}
\item The function $s_f:[0,\infty)\to[0,1]$ is increasing and subadditive, and $\lim_{\delta\to\infty}s_f(\delta)=1$.
\item The function $s_f$ is continuous on $[0,\infty)$.
\item The function $f$ is nonsingular if and only if $s_f(0)=0$.
\end{enumerate}
\end{lemma}

\begin{proof}
Let $\nu$ be the push-forward of Lebesgue measure on $[0,1]$ under $f$, that is,
\[
 \nu(A):=\mu(f^{-1}(A)).
\]
The monotonicity of $s_f$ is immediate. If $\mu(A)\leq\delta_1+\delta_2$, the non-atomicity of Lebesgue measure permits a decomposition $A=A_1\mathbin{\dot\cup}A_2$ with $\mu(A_i)\leq\delta_i$. Hence
\[
 \nu(A)\leq s_f(\delta_1)+s_f(\delta_2),
\]
and taking the supremum proves subadditivity. Furthermore,
\[
 \lim_{n\to\infty}\nu((-n,n))=\nu(\rr)=1,
\]
which implies $\lim_{\delta\to\infty}s_f(\delta)=1$.

For continuity, write the Lebesgue decomposition of $\nu$ as
\[
 \nu=\nu^{\mathrm{ac}}+\nu^{\mathrm{s}}
\]
with respect to Lebesgue measure. If $S$ is a null set supporting $\nu^{\mathrm{s}}$, then
\[
 s_f(0)=\nu^{\mathrm{s}}(\rr)
\]
and, for every $\delta\geq0$,
\begin{equation}\label{eq:stretch-decomp}
 s_f(\delta)=\nu^{\mathrm{s}}(\rr)+g(\delta),\qquad
 g(\delta):=\sup_{\mu(A)\leq\delta}\nu^{\mathrm{ac}}(A).
\end{equation}
Indeed, the upper bound follows from $\nu^{\mathrm{s}}(A)\leq\nu^{\mathrm{s}}(\rr)$, while the reverse bound is obtained by replacing $A$ with $A\cup S$. Since $\nu^{\mathrm{ac}}$ is a finite absolutely continuous measure, it is uniformly absolutely continuous: $g(t)\to0$ as $t\downarrow0$. The function $g$ is increasing and subadditive, and therefore, for $0\leq x\leq y$,
\[
 0\leq g(y)-g(x)\leq g(y-x).
\]
Thus $g$, and hence $s_f$, is continuous. Finally, $s_f(0)=0$ is precisely the statement that the preimage of every null set is null.
\end{proof}

The key point of Lemma~\ref{lem:stretch} is that for a nonsingular measurable function the preimages of sets of small measure are uniformly small: $\mu(f^{-1}(A))\leq s_f(\mu(A))$ and $s_f(t)\to0$ as $t\downarrow0$.

\subsection{Stretch function for polynomials}
To give an explicit algorithm that computes a transcendentally normal number, the general continuity of stretch functions is not precise enough. We instead need an a priori bound on $s_p$ for polynomials $p\in\zz[X]$.

\begin{proposition}\label{prop:poly}
Let $p\in\zz[X]$ be a non-constant polynomial of degree $d\geq1$. Then
\begin{equation}\label{eq:prop}
 s_p(\delta)\leq K_d\delta^{1/d},\qquad K_d:=4d.
\end{equation}
\end{proposition}

Before we turn to the proof, some remarks are in order.
\begin{enumerate}
\item The exponent $1/d$ is optimal, as can be seen by considering the monomial $x^d$.
\item The fact that the constant depends only on $d$ crucially uses $p\in\zz[X]$.
\item The explicit constant $K_d=4d$ is not optimal. Its precise value has only a minimal impact on the algorithm in Section~\ref{sec:alg}.
\end{enumerate}

Proposition~\ref{prop:poly} is based on the following form of P\'{o}lya's inequality.
\begin{lemma}\label{lem:polya}
Let $q\in\zz[X]+\rr$ be a polynomial of degree $d>0$ with integer coefficients except for a possible real constant term. Then
\begin{equation}\label{eq:polya}
 \mu\{x\in[0,1]:|q(x)|\leq\varepsilon^d\}\leq4\varepsilon.
\end{equation}
\end{lemma}

Lemma~\ref{lem:polya} is a standard consequence of P\'{o}lya's inequality \cite{Pol28}, or of Remez' inequality \cite{Rem36}. In Appendix~\ref{sec:polya} we give an elementary proof of a weaker bound which would also suffice for the computability arguments.

\begin{proof}[Proof of Proposition~\ref{prop:poly}]
The assertion is immediate for $d=1$. Suppose $d\geq2$, let $A\subset\rr$ be measurable with $\mu(A)\leq\delta$, and set
\[
 E:=p^{-1}(A)\cap[0,1],\qquad m:=\mu(E).
\]
If $m=0$ there is nothing to prove. The one-dimensional area formula and the fact that every value has at most $d$ preimages give
\begin{equation}\label{eq:area-poly}
 \int_E|p'(x)|\,dx
 =\int_A\#\{x\in E:p(x)=y\}\,dy
 \leq d\delta.
\end{equation}
Consequently, the set
\[
 E_0:=\left\{x\in E:|p'(x)|\leq\frac{2d\delta}{m}\right\}
\]
has measure at least $m/2$. Since $p'\in\zz[X]$ has degree $d-1$, Lemma~\ref{lem:polya} yields
\[
 \frac{m}{2}\leq4\left(\frac{2d\delta}{m}\right)^{1/(d-1)}.
\]
Hence
\[
 m^d\leq2d\,8^{d-1}\delta\leq(4d)^d\delta.
\]
Taking the supremum over $A$ proves \eqref{eq:prop}.
\end{proof}

\section{Construction of $\mathcal{F}$-normal Numbers }~\label{sec:const}

In this section we give an (abstract) explicit construction of $\mathcal{F}$-normal numbers for a family $\mathcal{F}$ of continuous nonsingular functions. To this end, we first discuss Sierpiński's method in the simpler setting of normal numbers and extend the technique in a second step.

\subsection{Sierpiński's cover method for normal numbers}

Sierpiński's method is based on the following elegant idea. Suppose we can find an open cover $U$ of the set of non-normal numbers $\mathcal{N}^c \cap [0,1]$ such that $W := [0,1] \setminus U \neq \emptyset.$ Then, $W$ is a closed nonempty set and in particular 
$$ \xi_W := \inf_{x \in W} x \in W \subset \mathcal{N},$$
i.e. $\xi_W$ is a concrete normal number! To construct an explicit normal, one "only" needs to find a suitable open cover $U$ - which is the core of Sierpiński's seminal work \cite{Sir17}. 

Let us recall Sierpiński's explicit open cover(s) (see also \cite{BF02}). We fix some $r \in (0,1]$ controlling the size of the cover and set
\begin{equation}\label{eq:sircover}
    U(r) := \bigcup_{b = 2}^{\infty} \bigcup_{m = 1}^{\infty}  \bigcup_{n = n_{m,b}(r)}^{\infty}  \bigcup_{d = 0}^{b-1} U_{b,m,n,d}
\end{equation}
with some union of intervals $U_{b,m,n,d}$ defined below. The integers $b  \geq 2$ correspond to the various bases, $m$ controls the fluctuation of the frequency of the digit $d$ from the "mean value" $\frac{1}{b}$ in the $b$-ary representation, the integer $n$ represents the length of the fractional representation under consideration and last but not least the integers $0 \leq d \leq b-1$ form the set of possible digits in a $b$-ary expansion of some number $x$. We introduce for any $d \in \zz$ the counting function $C_d : \zz^M \to \nn_0$
$$ C_d(q_1, \dots, q_M) := | \{ j = 1, \ldots M \, | \, q_j = d \}|. $$
The open set $U_{b,m,n,d}$ is then defined as 
\begin{equation}\label{eq:sirint}
    U_{b,m,n,d} :=\bigcup_{\substack{(q_1,\ldots, q_n) \in \{0,\ldots,b-1\}^n \\ |C_d(q_1, \dots, q_n) /n - 1/b| \geq 1/m }  } \left(\sum_{j=1}^{n} \frac{q_j}{b^j} \right)  + (-b^{-n}, 2 b^{-n}),
\end{equation}
where the addition is understood as translation of the interval $(-b^{-n}, 2 b^{-n})$ by the constant $\sum_{j=1}^{n} \frac{q_j}{b^j}$ to the right. Recall that a number is $b$-normal if and only if it is simply $b^k$-normal for every $k\geq1$ (see, for example, \cite{Knu97,Mad18}). Since the union in \eqref{eq:sircover} ranges over every integer base, it is therefore enough to cover failures of single-digit frequencies. In particular, the union $\bigcup_{n = n_{m,b}(r)}^{\infty}  U_{b,m,n,d}$ covers all real numbers $x$ for which the asymptotic frequency of the digit $d$ in its $b$-ary representation deviates more than $1/m$ from the "normal" value $1/b$. The key point is that this is true for any integer $n_{m,b}$. That is, one is free to choose $n_{m,b}$ as big as needed to control the measure of the whole open cover. In fact, the explicit choice 
\begin{equation}
    n_{m,b}(r) = \left\lfloor \frac{24m^6b^2}{r} \right\rfloor +2
\end{equation}
is good enough:
\begin{proposition}[\cite{BF02, Sir17}]\label{prop:sircov}
   Let $r \in (0,1]$ and $U(r)$ be the open set from \eqref{eq:sircover}. Then, $[0,1] \cap \mathcal{N}^c \subset U(r)$ and $\mu(U(r)) < r.$ In particular, the Sierpiński numbers
   \begin{equation}\label{eq:sirnum}
       \xi(r) := \inf \{ x \in [0,1] \, | \,x \notin U(r) \}
   \end{equation}
   exist and are normal.
\end{proposition}

A proof can be found in \cite{Sir17}. Our analysis in Appendix~\ref{sec:app1} essentially provides a simplified derivation of Proposition~\ref{prop:sircov}.

\subsection{Extending Sierpiński's approach to $\mathcal{F}$-normal numbers}

Based on Sierpiński's strategy, we give now a construction for an  $\mathcal{F}$-normal number. We will only consider countable families $\mathcal{F}$ of nonsingular continuous functions. Compared to the setting of Proposition~\ref{prop:full}, we need here to restrict to continuous functions in order to guarantee that preimages of open sets remain open. 

Let $r \in (0,1]$. We first extend the open cover $U(r)$ from \eqref{eq:sircover} to the whole real line:
\begin{equation}
    U^{+}(r) := \bigcup_{z \in \zz} z + U\left( \frac{r}{2^{|z| +2}}\right)
\end{equation}
Note that $\mu(U^{+}(r)) < r$ and $\mathcal{N}^c \subset  U^{+}(r). $ Let $f_1,f_2,\ldots$ be an enumeration of the countable set $\mathcal{F}$. Since each $f_k$ is nonsingular, Lemma~\ref{lem:stretch} allows us to make this choice monotonically. For $0<u<1$, set
\begin{equation}\label{eq:eta-choice}
 \eta_{f_k}(u):=\frac12\sup\{t\in(0,1]:s_{f_k}(t)<u\}.
\end{equation}
Then $\eta_{f_k}(u)>0$, the map $u\mapsto\eta_{f_k}(u)$ is increasing, and $s_{f_k}(\eta_{f_k}(u))<u$.
We define
\begin{equation}\label{eq:coverF}
 U_{\mathcal{F}}(r):=\bigcup_{k=1}^{\infty}f_k^{-1}\bigl(U^+(\eta_{f_k}(2^{-k}r))\bigr).
\end{equation}
By construction, $U_{\mathcal{F}}(r)$ is open and contains every $x\in[0,1]$ which is not $\mathcal{F}$-normal. Moreover,
\begin{align*}
 \mu(U_{\mathcal{F}}(r))
 &\leq\sum_{k=1}^{\infty}s_{f_k}\bigl(\mu(U^+(\eta_{f_k}(2^{-k}r)))\bigr)\\
 &<\sum_{k=1}^{\infty}2^{-k}r=r.
\end{align*}
Let us define the generalized Sierpiński numbers
\begin{equation}\label{eq:sirnumF}
    \xi_{\mathcal{F}}(r) := \inf \{ x \in [0,1] \, | \,x \notin U_{\mathcal{F}}(r) \},
\end{equation}
which are $\mathcal{F}$-normal and well-defined for $r \in (0,1]$ by the above. We summarize our findings in the following 
\begin{theorem}\label{thm:cover}
   Let $\mathcal{F}$ be a countable family of continuous and nonsingular functions. For any $r \in (0,1]$, the generalized Sierpiński number $\xi_{\mathcal{F}}(r) \in [0,1)$  is $\mathcal{F}$-normal.
\end{theorem}

Since the definition of $ \xi_{\mathcal{F}}(r)$ involves the inverses of $f_k$ and its stretch function $s_{f_k}$ it is considerably more involved and less explicit than the construction of the original Sierpiński numbers \eqref{eq:sirnum}. As our setting is rather general, one cannot expect an equally explicit procedure.

In the case of non-constant polynomials with integer coefficients, one can be more specific. Since integer translations preserve absolute normality, a number is t-normal if and only if it is $\mathcal{P}$-normal, where $\mathcal{P}$ denotes the family of polynomials
$$ \mathcal{P} := \{ p \in \zz[X] \setminus \{0 \} \, | \, p(0) = 0 \}. $$
In a first step, we want to fix our enumeration of the countable set $\mathcal{P}$. Let us introduce the norm 
$$ \left\| \sum_{k = 1}^d a_k x^k \right \|_{\mathcal{P}} := \sum_{k =1}^{d} 2^k | a_k|. $$
It is clear that for any number $K \in \nn$ the corresponding ball $B_K := \{ p \, | \, \| p \|_{\mathcal{P}} \leq K \}$ is a finite set. This allows us to define a total order $\preceq$ on $\mathcal{P}$ as follows:
\begin{align}\label{eq:order}
    \nonumber p = \sum_{k = 1}^{d_p} a_k x^k  &\preceq q = \sum_{k = 1}^{d_q} b_k x^k  \Longleftrightarrow  \| p \|_{\mathcal{P}} < \| q \|_{\mathcal{P}} \quad \text{ or } \\\quad &\| p \|_{\mathcal{P}} = \| q \|_{\mathcal{P}} \, \text{ and } a_{k_0} > b_{k_0} \text{ with } k_0 = \min\{ k \, | \, a_k \neq b_k\}
\end{align}
This order specifies a unique enumeration starting from the minimal element $p_1 = x$ and choosing always the well-defined successor. Moreover, in view of Proposition~\ref{prop:poly}, we may replace the implicit inverses of stretch functions, that is,
\begin{equation}\label{eq:coverP}
    U_{\mathcal{P}}(r) := \bigcup_{k =1}^{\infty} p^{-1}_k\left(U^{+}\left(\left(\frac{2^{-k}r}{K_{d_{p_k}}}\right)^{d_{p_k}}\right)\right),
\end{equation}
where $d_p$ denotes the degree of the polynomial and $K_d$ is the constant from Proposition~\ref{prop:poly}. For every $k$,
\[
 \mu\left([0,1]\cap p_k^{-1}\left(U^+\left(\left(\frac{2^{-k}r}{K_{d_{p_k}}}\right)^{d_{p_k}}\right)\right)\right)<2^{-k}r,
\]
and the same set contains every $x\in[0,1]$ for which $p_k(x)$ is non-normal. Hence $U_{\mathcal P}(r)$ covers the complement of the t-normal numbers in $[0,1]$ and has measure less than $r$.
Note that this definition does not coincide with \eqref{eq:coverF} for $\mathcal{F} = \mathcal{P}.$ Accordingly, one defines the Sierpiński numbers 
\begin{equation}\label{eq:sirnumP}
    \xi_{\mathcal{P}}(r) := \inf \{ x \in [0,1] \, | \,x \notin U_{\mathcal{P}}(r) \},
\end{equation}

Apart from the algebraic inverses in $\eqref{eq:coverP}$, this construction is as explicit as Sierpiński's original one.
\subsection{Properties of Sierpiński numbers}
In the following proposition, we collect two elementary properties of the Sierpiński numbers $\xi_{\mathcal{F}}(r)$.
\begin{proposition}\label{prop:sirprop}
Let $\mathcal{F}$ be a countable family of continuous nonsingular functions. Then:
\begin{enumerate}
\item The map $(0,1]\ni r\mapsto\xi_{\mathcal{F}}(r)$ is increasing and $\lim_{r\downarrow0}\xi_{\mathcal{F}}(r)=0$.
\item If $0\notin\mathcal{N}_{\mathcal{F}}$, then the family $(\xi_{\mathcal{F}}(r))_{r\in(0,1]}$ is infinite.
\end{enumerate}
\end{proposition}

\begin{proof}
By construction, $U_{\mathcal{F}}(r)\subset U_{\mathcal{F}}(r')$ whenever $0<r\leq r'\leq1$, and hence $\xi_{\mathcal{F}}(r)\leq\xi_{\mathcal{F}}(r')$. Moreover,
\[
 [0,\xi_{\mathcal{F}}(r))\subset U_{\mathcal{F}}(r),
\]
so that
\begin{equation}\label{eq:sirbound}
 \xi_{\mathcal{F}}(r)\leq\mu(U_{\mathcal{F}}(r))<r.
\end{equation}
This proves the first assertion. If $0$ is not $\mathcal{F}$-normal, then $0\in U_{\mathcal{F}}(r)$; since the latter set is open, $\xi_{\mathcal{F}}(r)>0$. Define
\[
 \xi_1:=\xi_{\mathcal{F}}(1/2),\qquad
 \xi_{n+1}:=\xi_{\mathcal{F}}(\xi_n).
\]
By \eqref{eq:sirbound}, $0<\xi_{n+1}<\xi_n$, which proves infinitude.
\end{proof}

Let us turn to the computational perspective, where we restrict ourselves to $\mathcal{F}=\mathcal{P}$. The algorithm in Section~\ref{sec:alg} does not compute a Sierpiński number, but another t-normal number. For rational parameters, the Sierpiński numbers are nevertheless computably enumerable.

\begin{definition}[\cite{BF02}]\label{def:comenu}
A real number $r$ is called computably enumerable if there is an algorithm producing rational numbers $q_n$ such that $q_n\leq q_{n+1}$ and $q_n\to r$.
\end{definition}

\begin{proposition}\label{prop:comenu}
The Sierpiński number $\xi_{\mathcal{P}}(r)$ is computably enumerable for every rational $r\in(0,1]$.
\end{proposition}

\begin{proof}
For rational $r$, the open set $U_{\mathcal{P}}(r)$ admits a computable enumeration
\[
 U_{\mathcal{P}}(r)=\bigcup_{k=1}^{\infty}p_{n(k)}^{-1}(I_k),
\]
where the $I_k$ have rational endpoints. Let
\[
 U_K:=\bigcup_{k=1}^{K}p_{n(k)}^{-1}(I_k),\qquad
 \xi_K:=\min([0,1]\setminus U_K).
\]
The complement is nonempty because $\mu(U_K)<1$. The sequence $\xi_K$ is increasing and bounded above by $\xi_{\mathcal{P}}(r)$. If $\xi_*:=\lim_K\xi_K$, then, for every fixed $K$, one has $\xi_M\in[0,1]\setminus U_K$ for $M\geq K$. Since this complement is closed, $\xi_*\notin U_K$ for every $K$, and therefore $\xi_*\notin U_{\mathcal{P}}(r)$. Hence $\xi_*=\xi_{\mathcal{P}}(r)$.

Each $U_K$ is a finite union of intervals with computable algebraic endpoints, obtained by isolating the real roots of finitely many polynomials with rational coefficients. Thus $\xi_K$ is computable. Choose a rational $q_K$ with $\xi_K-2^{-K}<q_K\leq\xi_K$ and set $\widehat q_K:=\max_{j\leq K}q_j$. Then $(\widehat q_K)$ is increasing and converges to $\xi_{\mathcal{P}}(r)$.
\end{proof}

\section{Construction of $\F^\lil$-normal numbers}\label{sec:lil}

We turn to the more sophisticated covers for LIL-normal numbers. As in Section~\ref{sec:const}, the passage from one exceptional set to a countable family of functions is handled by stretch functions.

\begin{corollary}\label{cor:FLIL}
Let $(U^\lil(r))_{r\in(0,1]}$ be a family of open sets such that
\begin{equation}\label{eq:coverLILprop}
 (\mathcal N^\lil)^c\cap[0,1]\subset U^\lil(r),\qquad \mu(U^\lil(r))<r.
\end{equation}
Given a countable family $\F=(f_k)_{k\in\nn}$ of continuous nonsingular functions, set
\begin{align}\label{eq:coverFLIL}
 U_{\F}^\lil(r)&:=\bigcup_{k=1}^{\infty}f_k^{-1}\bigl(U_+^\lil(\eta_{f_k}(2^{-k}r))\bigr),\\
 U_+^\lil(s)&:=\bigcup_{z\in\zz}\left(z+U^\lil\left(\frac{s}{2^{|z|+2}}\right)\right),\nonumber
\end{align}
where $\eta_{f_k}$ is chosen as in \eqref{eq:eta-choice}. Then
\begin{equation}\label{eq:sirnumFLIL}
 \xi_{\F}^\lil(r):=\inf\{x\in[0,1]:x\notin U_{\F}^\lil(r)\}
\end{equation}
is $\F^\lil$-normal. In the polynomial case one may use the explicit cover
\begin{equation}\label{eq:coverPLIL}
 U_{\mathcal P}^\lil(r):=\bigcup_{k=1}^{\infty}p_k^{-1}\left(
 U_+^\lil\left(\left(\frac{2^{-k}r}{K_{d_{p_k}}}\right)^{d_{p_k}}\right)\right).
\end{equation}
The corresponding numbers
\begin{equation}\label{eq:sirnumPLIL}
 \xi_{\mathcal P}^\lil(r):=\inf\{x\in[0,1]:x\notin U_{\mathcal P}^\lil(r)\}
\end{equation}
are transcendentally LIL-normal.
\end{corollary}

\begin{proof}
Integer translations preserve the fractional parts of $b^{j-1}x$ for every integer base $b$. Thus $U_+^\lil(s)$ covers all non-LIL-normal real numbers and has measure less than $s$. The stretch-function estimate gives
\[
 \mu(U_{\F}^\lil(r))<\sum_{k=1}^{\infty}2^{-k}r=r.
\]
The polynomial assertion follows from Proposition~\ref{prop:poly}; here the power $d_{p_k}$ is essential.
\end{proof}

We can further reduce the construction to a fixed base.
\begin{lemma}\label{lem:LILb}
Suppose that, for every integer $b\geq2$ and $r\in(0,1]$, the set $U_b^\lil(r)$ is open, has measure less than $r$, and every $x\in[0,1]\setminus U_b^\lil(r)$ satisfies the LIL in base $b$. Then
\begin{equation}\label{eq:all-bases-cover}
 U^\lil(r):=\bigcup_{b=2}^{\infty}U_b^\lil(2^{-b}r)
\end{equation}
meets the assumptions of Corollary~\ref{cor:FLIL}.
\end{lemma}

\begin{proof}
The set is open, it covers $(\mathcal N^\lil)^c\cap[0,1]$, and
\[
 \mu(U^\lil(r))<r\sum_{b=2}^{\infty}2^{-b}=\frac r2<r.
\]
\end{proof}

\subsection{Outline for the construction of $U_b^\lil(r)$}
Fix an integer base $b\geq2$. For $0\leq a<a'\leq1$, write
\begin{equation}\label{eq:dispar1}
 D_N^{a,a'}(x):=\left|\frac1N\sum_{j=1}^{N}\mathbbm{1}_{[a,a')}\bigl(\{b^{j-1}x\}\bigr)-(a'-a)\right|.
\end{equation}
Fukuyama's proof of the discrepancy LIL consists of two ingredients. First, for fixed $a,a'$ one has
\begin{equation}\label{eq:lilres}
 \limsup_{N\to\infty}\frac{\sqrt N D_N^{a,a'}(x)}{\sqrt{2\log\log N}}=\sigma_{a,a'}
 \quad\text{for a.e. }x,
\end{equation}
where $\sigma_{a,a'}$ is explicit. Secondly, the contribution of intervals whose length tends to zero is negligible:
\begin{equation}\label{eq:discr}
 \lim_{L\to\infty}\limsup_{N\to\infty}
 \sup_{|a-a'|\leq2^{-L}}
 \frac{\sqrt N D_N^{a,a'}(x)}{\sqrt{2\log\log N}}=0
 \quad\text{for a.e. }x.
\end{equation}
We make both ingredients effective. For $b$-adic endpoints, the sequence of indicators is finitely dependent and can be split into independent blocks. The second ingredient is obtained from the Fukuyama--Philipp discretization argument.

The construction will prove the following theorem.
\begin{theorem}\label{thm:ulilb}
Let $r\in(0,1]$, and let $\mathcal L_b^{a,a'}(s)$ and $\mathcal D_b(s)$ be the open sets defined in \eqref{eq:L} and \eqref{eq:D}. Set
\begin{equation}\label{eq:covUb}
 U_b^\lil(r):=\mathcal D_b\left(\frac r2\right)\cup
 \bigcup_{L=1}^{\infty}\ \bigcup_{0\leq k<k'<b^L}
 \mathcal L_b^{k/b^L,k'/b^L}\left(2^{-(L+1)}b^{-2L}r\right).
\end{equation}
Then $U_b^\lil(r)$ is open, $\mu(U_b^\lil(r))<r$, and every $x\in[0,1]\setminus U_b^\lil(r)$ satisfies
\begin{equation}\label{eq:full-LIL-conclusion}
 \limsup_{N\to\infty}\frac{\sqrt N D_N^b(x)}{\sqrt{\log\log N}}=L_b,
\end{equation}
where
\[
 L_b=\begin{cases}
 \sqrt{84/81},&b=2,\\
 \sqrt{(b+1)/(2(b-1))},&b\text{ odd},\\
 \sqrt{b(b+1)(b-2)/(2(b-1)^3)},&b\geq4\text{ even}.
 \end{cases}
\]
\end{theorem}

\subsection{Effective bounds for the standard LIL}
For $N\geq3$, set
\[
 \phi_N:=\sqrt{2N\log\log N}.
\]
The following proposition gives the effective form of the standard LIL estimates needed below.

\begin{proposition}\label{prop:RW}
Let $Y_1,Y_2,\ldots$ be i.i.d. bounded random variables, let $\mu=\ee Y_1$, $\sigma^2=\operatorname{Var}(Y_1)>0$, and set $S_N=\sum_{j=1}^{N}Y_j$. For $\delta>0$, define
\begin{equation}\label{eq:VLILup}
 V_{\delta,N_0}^{\mathrm{up}}:=\bigcup_{N\geq N_0}
 \{|S_N-\mu N|>(1+\delta)\sigma\phi_N\}.
\end{equation}
For $0<\delta<1$, define
\begin{equation}\label{eq:VLILlow}
 V_{\delta,N_0}^{\mathrm{low}}:=\bigcup_{N\geq N_0}
 \bigcap_{k=2}^{N}\{|S_{N^k}-\mu N^k|<(1-\delta)\sigma\phi_{N^k}\}.
\end{equation}
There are explicit tail functions $R_\delta^{\mathrm{up}}(N_0)$ and $R_\delta^{\mathrm{low}}(N_0)$, depending only on $\delta$, $\sigma$, a bound for $|Y_1-\mu|$, and the third absolute centered moment, such that
\begin{equation}\label{eq:RW-tail-bounds}
 \pp(V_{\delta,N_0}^{\mathrm{up}})\leq R_\delta^{\mathrm{up}}(N_0),\qquad
 \pp(V_{\delta,N_0}^{\mathrm{low}})\leq R_\delta^{\mathrm{low}}(N_0),
\end{equation}
and both right-hand sides tend to zero as $N_0\to\infty$. More precisely, for suitable constants $c_i=c_i(\delta)>0$,
\begin{align}
 R_\delta^{\mathrm{up}}(N_0)&\leq c_1\sum_{j\geq c_2\log N_0}j^{-1-c_3},\label{eq:Rup-form}\\
 R_\delta^{\mathrm{low}}(N_0)&\leq\sum_{N\geq N_0}
 \exp\left(-c_1\frac{N^{c_3}}{(\log N)^{c_2}}\right).
 \label{eq:Rlow-form}
\end{align}
If the law of $Y_1$ has computable finite support, these bounds and a valid starting index are computable uniformly from the law and $\delta$.
\end{proposition}

\begin{proof}
Put $X_j:=Y_j-\mu$, $T_n:=\sum_{j=1}^{n}X_j$, and choose $B$ such that $|X_1|\leq B$. We use two standard inequalities (see, for example, \cite{Pet75}). First, the maximal Bernstein inequality gives, for $t>0$,
\begin{equation}\label{eq:max-bernstein}
 \pp\left(\max_{m\leq n}|T_m|\geq t\right)
 \leq2\exp\left(-\frac{t^2}{2(n\sigma^2+Bt/3)}\right).
\end{equation}
Secondly, the Berry--Esseen theorem gives, with $\rho_3:=\ee|X_1|^3$ and the standard normal distribution function $\Phi$,
\begin{equation}\label{eq:berry-esseen}
 \sup_x\left|\pp\left(\frac{T_m}{\sigma\sqrt m}\leq x\right)-\Phi(x)\right|
 \leq\frac{C_{\mathrm{BE}}\rho_3}{\sigma^3\sqrt m}.
\end{equation}
We also use the elementary lower Mills bound
\begin{equation}\label{eq:mills-lower}
 1-\Phi(x)\geq c_0\frac{e^{-x^2/2}}{1+x},\qquad x\geq0,
\end{equation}
with a universal $c_0>0$.

For the upper bound, choose $\theta>1$ and $\varepsilon>0$ so small that
\[
 c_*:=\frac{(1+\delta)^2}{\theta(1+\varepsilon)}>1,
\]
and set $n_j=\lfloor\theta^j\rfloor$. For all sufficiently large $j$, monotonicity of $\phi_n$, \eqref{eq:max-bernstein}, and $n_{j+1}\leq\theta n_j+1$ imply
\begin{align*}
 &\pp\left(\max_{n_j\leq m<n_{j+1}}|T_m|>(1+\delta)\sigma\phi_m\right)\\
 &\qquad\leq\pp\left(\max_{m\leq n_{j+1}}|T_m|>(1+\delta)\sigma\phi_{n_j}\right)
 \leq C j^{-c_*}.
\end{align*}
Here the term $Bt/3$ in \eqref{eq:max-bernstein} is absorbed into the factor $1+\varepsilon$ for large $j$, and $\log\log n_j=\log j+O(1)$. Summation over the blocks proves \eqref{eq:Rup-form}.

For the lower bound, let
\[
 B_N:=\bigcap_{k=2}^{N}\{|T_{N^k}|<(1-\delta)\sigma\phi_{N^k}\}.
\]
For $2\leq k\leq N-1$, define the independent increments
\[
 Z_{N,k}:=T_{N^{k+1}}-T_{N^k},\qquad m_{N,k}:=N^{k+1}-N^k.
\]
On $B_N$,
\begin{equation}\label{eq:block-small}
 |Z_{N,k}|\leq(1-\delta)\sigma(\phi_{N^{k+1}}+\phi_{N^k}).
\end{equation}
Uniformly for $2\leq k\leq N-1$, and for all sufficiently large $N$,
\begin{equation}\label{eq:xNk-bound}
 \frac{(1-\delta)(\phi_{N^{k+1}}+\phi_{N^k})}{\sqrt{m_{N,k}}}
 \leq(1-\delta/2)\sqrt{2\log((k+1)\log N)}.
\end{equation}
Set $\beta:=(1-\delta/2)^2<1$. Applying \eqref{eq:berry-esseen} and \eqref{eq:mills-lower} to $Z_{N,k}$, and using $m_{N,k}\geq N^2(N-1)$, yields, uniformly in $2\leq k\leq N-1$,
\begin{equation}\label{eq:block-exceed}
 \pp\left(|Z_{N,k}|>(1-\delta)\sigma(\phi_{N^{k+1}}+\phi_{N^k})\right)
 \geq c\frac{((k+1)\log N)^{-\beta}}
 {\sqrt{\log((k+1)\log N)}}
\end{equation}
for large $N$. Indeed, the Berry--Esseen error is $O(m_{N,k}^{-1/2})$, which is uniformly negligible compared with the displayed Gaussian lower bound. Independence and $1-x\leq e^{-x}$ now give
\begin{align*}
 \pp(B_N)
 &\leq\exp\left(-c\sum_{k=2}^{N-1}
 \frac{((k+1)\log N)^{-\beta}}
 {\sqrt{\log((k+1)\log N)}}\right).
\end{align*}
Restricting the sum to $\lceil N/2\rceil\leq k\leq N-1$ gives at least $N/2-1$ terms, each bounded below by a constant multiple of
$N^{-\beta}(\log N)^{-\beta-1/2}$. Thus, after increasing the starting index,
\[
 \pp(B_N)\leq\exp\left(-c'\frac{N^{1-\beta}}
 {(\log N)^{\beta+1/2}}\right).
\]
The union bound over $N\geq N_0$ proves \eqref{eq:Rlow-form}. All constants arise from explicit inequalities. For a finite computable law, the moments are computable and one can search for a starting index at which the finitely many displayed elementary inequalities hold.
\end{proof}

\subsection{A cover for $b$-adic intervals}\label{sec:bcover}
Fix two $b$-adic numbers
\[
 a=\frac{q_a}{b^L}<a'=\frac{q_{a'}}{b^L},\qquad 0\leq q_a<q_{a'}<b^L.
\]
Define
\begin{equation}\label{eq:block-indicator}
 w_j(x):=\left\lfloor b^L\{b^{j-1}x\}\right\rfloor,
 \qquad X_j(x):=\mathbbm{1}_{[q_a,q_{a'})}(w_j(x)).
\end{equation}
Then
\begin{equation}\label{eq:dispar2}
 D_N^{a,a'}(x)=\left|\frac1N\sum_{j=1}^{N}X_j(x)-(a'-a)\right|.
\end{equation}
Under Lebesgue measure, $(X_j)$ is stationary, and every subcollection $(X_j)_{j\in J}$ with pairwise distances at least $L$ is independent.

For $M\in\nn$, let
\begin{align*}
 J_k&:=\{(k-1)L(M+1)+1,\ldots,(k-1)L(M+1)+LM\},\\
 \widehat J_k&:=\{(k-1)L(M+1)+LM+1,\ldots,kL(M+1)\},
\end{align*}
and set
\begin{equation}\label{eq:YZ-blocks}
 Y_k^M:=\sum_{j\in J_k}X_j,
 \qquad Z_k^M:=\sum_{j\in\widehat J_k}X_j.
\end{equation}
Both $(Y_k^M)_k$ and $(Z_k^M)_k$ are i.i.d.; the law of $Z_k^M$ is the same as that of $Y_k^1$.

\begin{lemma}\label{lem:varY}
Let
\[
 \sigma_M^2:=\frac1{ML}\ee\bigl[(Y_1^M-ML(a'-a))^2\bigr].
\]
Then
\begin{equation}\label{eq:varM}
 \lim_{M\to\infty}\sigma_M^2=\sigma_{a,a'}^2,
\end{equation}
where
\begin{equation}\label{eq:sigma-aa}
 \sigma_{a,a'}^2:=\int_0^1\psi_{a,a'}(x)^2\,dx
 +2\sum_{j=1}^{L-1}\int_0^1\psi_{a,a'}(x)\psi_{a,a'}(b^jx)\,dx
\end{equation}
and $\psi_{a,a'}(x)=\mathbbm{1}_{[a,a')}(\{x\})-(a'-a)$.
\end{lemma}

\begin{proof}
Expanding the square and using the independence of $X_k$ and $X_{k'}$ for $|k-k'|\geq L$ gives
\begin{align*}
 \ee[(Y_1^M-ML(a'-a))^2]
 &=LM\int_0^1\psi_{a,a'}^2(x)\,dx\\
 &\quad+\sum_{j=1}^{L-1}(2LM-2j)
 \int_0^1\psi_{a,a'}(x)\psi_{a,a'}(b^jx)\,dx.
\end{align*}
Division by $LM$ proves the claim.
\end{proof}

We next define open exceptional sets. The rational cutoffs below avoid any decision problem caused by equality with a computable real threshold. The random variables $Y_1^M$ and $Z_1^M$ are non-constant whenever $0<a'-a<1$: the constant base-$b$ word $(b-1)(b-1)\cdots$ makes every summand vanish, whereas a word whose first $L$ digits represent an integer in $[q_a,q_{a'})$ makes at least one summand equal to one. Hence all variances below are strictly positive. For $M\geq2$ and $q\geq3$, choose computable rational numbers $u_{M,q}^+$, $u_{M,q}^-$, and $v_q^+$ such that
\begin{align}
 \left(1+\frac1{2M}\right)\sqrt{LM}\sigma_M\phi_q
 &<u_{M,q}^+<\left(1+\frac1M\right)\sqrt{LM}\sigma_M\phi_q,\label{eq:rational-up}\\
 \left(1-\frac1M\right)\sqrt{LM}\sigma_M\phi_q
 &<u_{M,q}^-<\left(1-\frac1{2M}\right)\sqrt{LM}\sigma_M\phi_q,\label{eq:rational-low}\\
 \frac32\sqrt L\sigma_1\phi_q&<v_q^+<2\sqrt L\sigma_1\phi_q.\label{eq:rational-z}
\end{align}
Let
\begin{align*}
 E_{M,q}^{Y,+}&:=\left\{\left|\sum_{j=1}^{q}Y_j^M-LMq(a'-a)\right|>u_{M,q}^+\right\},\\
 E_{M,q,k}^{Y,-}&:=\left\{\left|\sum_{j=1}^{q^k}Y_j^M-LMq^k(a'-a)\right|<u_{M,q^k}^-\right\},\\
 E_{M,q}^{Z,+}&:=\left\{\left|\sum_{j=1}^{q}Z_j^M-Lq(a'-a)\right|>v_q^+\right\}.
\end{align*}
These are finite unions of half-open intervals with rational endpoints. Put
\begin{align*}
 \rho^Y_{M,q}&:=2^{-q-2}b^{-(M+1)qL+1},\\
 \widehat\rho^Y_{M,q,k}&:=2^{-q-k-2}b^{-(M+1)q^kL+1},\\
 \rho^Z_{M,q}&:=2^{-q-2}b^{-(M+1)qL-L+1}.
\end{align*}
Define
\begin{align}
 \mathcal Y_{N,M,b}^{\mathrm{up}}(a,a')
 &:=\bigcup_{q\geq N}\left(E_{M,q}^{Y,+}+(-\rho^Y_{M,q},\rho^Y_{M,q})\right),\label{eq:Yup}\\
 \mathcal Y_{N,M,b}^{\mathrm{low}}(a,a')
 &:=\bigcup_{q\geq N}\ \bigcap_{k=2}^{q}
 \left(E_{M,q,k}^{Y,-}+(-\widehat\rho^Y_{M,q,k},\widehat\rho^Y_{M,q,k})\right),\label{eq:Ylow}\\
 \mathcal Z_{N,M,b}^{\mathrm{up}}(a,a')
 &:=\bigcup_{q\geq N}\left(E_{M,q}^{Z,+}+(-\rho^Z_{M,q},\rho^Z_{M,q})\right).
 \label{eq:Zup}
\end{align}
The additions are Minkowski sums. The first event at level $q$ is measurable with respect to the first $(M+1)qL-1$ base-$b$ digits, the event in \eqref{eq:Ylow} at $(q,k)$ with respect to the first $(M+1)q^kL-1$ digits, and the last event with respect to the first $(M+1)qL+L-1$ digits. Consequently, the total measure added by the neighborhoods in each of \eqref{eq:Yup}--\eqref{eq:Zup} is at most $2^{-N}$. Proposition~\ref{prop:RW} and \eqref{eq:rational-up}--\eqref{eq:rational-z} therefore imply
\begin{align}
 \mu(\mathcal Y_{N,M,b}^{\mathrm{up}}(a,a'))
 &\leq R_{Y^M,1/(2M)}^{\mathrm{up}}(N)+2^{-N},\label{eq:Yup-bound}\\
 \mu(\mathcal Y_{N,M,b}^{\mathrm{low}}(a,a'))
 &\leq R_{Y^M,1/(2M)}^{\mathrm{low}}(N)+2^{-N},\label{eq:Ylow-bound}\\
 \mu(\mathcal Z_{N,M,b}^{\mathrm{up}}(a,a'))
 &\leq R_{Z^M,1/2}^{\mathrm{up}}(N)+2^{-N}.
 \label{eq:Zup-bound}
\end{align}
Here the subscript indicates the increment law to which Proposition~\ref{prop:RW} is applied.

For $r\in(0,1)$ choose $N_M(r)$ so large that the sum of the three right-hand sides in \eqref{eq:Yup-bound}--\eqref{eq:Zup-bound} is smaller than $2^{-M}r$. If $r$ is computable, this choice is effective. Set
\begin{equation}\label{eq:L}
 \mathcal L_b^{a,a'}(r):=\bigcup_{M=2}^{\infty}
 \left(\mathcal Y_{N_M(r),M,b}^{\mathrm{up}}(a,a')
 \cup\mathcal Y_{N_M(r),M,b}^{\mathrm{low}}(a,a')
 \cup\mathcal Z_{N_M(r),M,b}^{\mathrm{up}}(a,a')\right).
\end{equation}

\begin{proposition}\label{prop:L}
The set $\mathcal L_b^{a,a'}(r)$ is open, has measure less than $r$, and every $x\in[0,1]\setminus\mathcal L_b^{a,a'}(r)$ satisfies
\begin{equation}\label{eq:LILa}
 \limsup_{N\to\infty}\frac{\sqrt N D_N^{a,a'}(x)}{\sqrt{2\log\log N}}=\sigma_{a,a'}.
\end{equation}
For computable $r$, the defining union and its tail bounds are computable.
\end{proposition}

\begin{proof}
Openness and the measure estimate follow from the construction. Fix $M\geq2$ and put
\[
 T_n(x):=\sum_{j=1}^{n}(X_j(x)-(a'-a)),\qquad n_q:=L(M+1)q.
\]
At the cycle endpoints,
\[
 T_{n_q}=\sum_{j=1}^{q}(Y_j^M-LM(a'-a))
 +\sum_{j=1}^{q}(Z_j^M-L(a'-a)).
\]
If $x$ lies outside the three exceptional sets associated with $M$, then, for all sufficiently large $q$,
\begin{align*}
 \left|\sum_{j=1}^{q}(Y_j^M-LM(a'-a))\right|
 &\leq\left(1+\frac1M\right)\sqrt{LM}\sigma_M\phi_q,\\
 \left|\sum_{j=1}^{q}(Z_j^M-L(a'-a))\right|
 &\leq2\sqrt L\sigma_1\phi_q.
\end{align*}
Moreover, for every sufficiently large integer $Q$ there is a $k\in\{2,\ldots,Q\}$ such that, with $q=Q^k$,
\[
 \left|\sum_{j=1}^{q}(Y_j^M-LM(a'-a))\right|
 \geq\left(1-\frac1M\right)\sqrt{LM}\sigma_M\phi_q.
\]
The difference between $T_n$ and the nearest cycle endpoint is bounded by $L(M+1)$, hence it is negligible after division by $\sqrt{2n\log\log n}$. Since
\[
 \frac{\phi_q}{\sqrt{2n_q\log\log n_q}}\longrightarrow\frac1{\sqrt{L(M+1)}},
\]
we obtain
Writing
\begin{align*}
 A_M&:=\sqrt{\frac M{M+1}}\left(1-\frac1M\right)\sigma_M
 -\frac{2}{\sqrt{M+1}}\sigma_1,\\
 B_M&:=\sqrt{\frac M{M+1}}\left(1+\frac1M\right)\sigma_M
 +\frac{2}{\sqrt{M+1}}\sigma_1,
\end{align*}
we obtain
\[
 A_M\leq\limsup_{N\to\infty}
 \frac{|T_N(x)|}{\sqrt{2N\log\log N}}\leq B_M.
\]
Letting $M\to\infty$ and using Lemma~\ref{lem:varY} proves \eqref{eq:LILa}. The splitting into independent subsequences is closely related to the classical arguments of Kac \cite{Kac46} and Maruyama \cite{Mar50}.
\end{proof}

\subsection{The Fukuyama--Philipp discretization argument}
We now make the small-interval estimate \eqref{eq:discr} effective. We first recall the deterministic decomposition used by Philipp. Fix $L\in\nn$ and $K\in\{0,\ldots,2^L-1\}$, and write $a(K,L)=K2^{-L}$. For $n$ large, let
\[
 H_n:=\left\lfloor\frac n2\right\rfloor+1
\]
and let $\epsilon=(\epsilon_{L+1},\ldots,\epsilon_{H_n})\in\{0,1\}^{H_n-L}$. Put
\[
 u_h:=\frac K{2^L}+\sum_{j=L+1}^{h}\frac{\epsilon_j}{2^j},\qquad L\leq h\leq H_n,
\]
with the empty sum interpreted as zero, and define the mean-zero step functions
\begin{equation}\label{eq:rho}
 \rho_{h,\epsilon}^{L,K}(x):=
 \begin{cases}
 \mathbbm{1}_{[u_h,u_{h+1})}(\{x\})-(u_{h+1}-u_h),&L\leq h<H_n,\\
 \mathbbm{1}_{[u_{H_n},u_{H_n}+2^{-H_n})}(\{x\})-2^{-H_n},&h=H_n.
 \end{cases}
\end{equation}
Thus the centering in the first line is exactly $2^{-(h+1)}\epsilon_{h+1}$. Set
\begin{equation}\label{eq:fukF}
 F_{M,N,h,\epsilon}^{L,K}(x):=
 \left|\sum_{j=M+1}^{M+N}\rho_{h,\epsilon}^{L,K}(b^jx)\right|.
\end{equation}
The following is the deterministic part of Philipp's argument \cite[Lemma~4]{Phi75}.

\begin{lemma}\label{lem:philipp-decomp}
Let $a(K,L)<a'\leq a(K,L)+2^{-L}$ and $2^n\leq N<2^{n+1}$. There are a bit string $\epsilon$ and integers $0\leq m_\ell<2^{n-\ell}$ such that
\begin{equation}\label{eq:fukbound}
 N D_N^{a(K,L),a'}(x)
 \leq\sum_{h=L}^{H_n}\left(
 F_{0,2^n,h,\epsilon}^{L,K}(x)
 +\sum_{\lceil n/3\rceil\leq\ell\leq n}
 F_{2^n+m_\ell2^\ell,2^{\ell-1},h,\epsilon}^{L,K}(x)
 \right)+2\sqrt N.
\end{equation}
\end{lemma}

For algorithmic purposes, we again use rational cutoffs. Fix the absolute constant $C_P$ occurring in Philipp's estimates. Choose rational numbers $\tau_{n,h}$ and $\widehat\tau_{n,h,\ell}$ satisfying
\begin{align*}
 C_P2^{-h/8}\phi_{2^n}<\tau_{n,h}<(C_P+1)2^{-h/8}\phi_{2^n},\\
 C_P2^{-h/8}2^{\ell-n-3}\phi_{2^n}<\widehat\tau_{n,h,\ell}
 <(C_P+1)2^{-h/8}2^{\ell-n-3}\phi_{2^n}.
\end{align*}
Define
\begin{align*}
 G^{L,K}(n,\epsilon,h)&:=\{x\in[0,1]:F_{0,2^n,h,\epsilon}^{L,K}(x)\geq\tau_{n,h}\},\\
 \widehat G^{L,K}(n,\epsilon,h,\ell,m)&:=
 \{x\in[0,1]:F_{2^n+m2^\ell,2^{\ell-1},h,\epsilon}^{L,K}(x)
 \geq\widehat\tau_{n,h,\ell}\},
\end{align*}
and let $E^{L,K}(n_0)$ be the union of these sets over
\[
 n\geq n_0,\quad L\leq h\leq H_n,\quad
 \epsilon\in\{0,1\}^{H_n-L},\quad
 \lceil n/3\rceil\leq\ell\leq n,\quad0\leq m<2^{n-\ell},
\]
where the $\widehat G$-terms are omitted in the first part of the union.

Each individual set is a finite union of half-open intervals with rational endpoints, uniformly computable from its indices. If $A$ is such a finite union and $\eta>0$ is rational, let $\mathcal O(A,\eta)$ denote the open set obtained after merging its components and enlarging each of the resulting $q$ intervals by $\eta/(4q)$ on both sides; by convention, $\mathcal O(\varnothing,\eta)=\varnothing$. Then
\begin{equation}\label{eq:open-enlargement}
 A\subset\mathcal O(A,\eta),\qquad
 \mu(\mathcal O(A,\eta))\leq\mu(A)+\eta.
\end{equation}
At every fixed level $n$, let $Q_n^{L,K}$ be the number of sets occurring above, and define $E_+^{L,K}(n_0)$ by replacing every level-$n$ set $A$ with
\[
 \mathcal O\left(A,\frac{2^{-n-1}}{Q_n^{L,K}}\right).
\]
This replacement is preferable to assigning a single radius before the number of interval components has been counted.

\begin{lemma}\label{lem:phi}
There are absolute effectively computable constants $C_P,C_\Phi>0$ and $n_\Phi\in\nn$ such that, for $n_0\geq\max\{n_\Phi,2L\}$,
\begin{equation}\label{eq:phi2}
 \mu(E^{L,K}(n_0))\leq\frac{C_\Phi}{n_0}.
\end{equation}
Moreover,
\begin{equation}\label{eq:phi1}
 \mu(E_+^{L,K}(n_0))\leq\frac{C_\Phi}{n_0}+2^{-n_0}.
\end{equation}
For $x\notin E_+^{L,K}(n_0)$,
\begin{equation}\label{eq:philipp-small-cell}
 \limsup_{N\to\infty}\sup_{a(K,L)<a'\leq a(K,L)+2^{-L}}
 \frac{\sqrt N D_N^{a(K,L),a'}(x)}{\sqrt{2\log\log N}}
 \leq C2^{-L/8}
\end{equation}
with an absolute constant $C$.
\end{lemma}

\begin{proof}
The estimate \eqref{eq:phi2} is the quantitative content of Philipp's Lemma~5 \cite{Phi75}; all constants in that proof are explicit and uniform for integer $b\geq2$. By \eqref{eq:open-enlargement}, the total increase of measure at level $n$ is at most $2^{-n-1}$. Summing over $n\geq n_0$ proves \eqref{eq:phi1}.

If $x\notin E_+^{L,K}(n_0)$, then none of the original exceptional events occurs. Inserting the rational upper bounds for the functions $F$ into \eqref{eq:fukbound}, summing the geometric factors $2^{-h/8}$ and $2^{\ell-n-3}$, and using $\phi_{2^n}\asymp\phi_N$, gives \eqref{eq:philipp-small-cell}. The term $2\sqrt N$ vanishes after normalization.
\end{proof}

For $r\in(0,1)$, choose $n_L(r)\geq\max\{n_\Phi,2L\}$ computably so that
\begin{equation}\label{eq:nL-choice}
 \frac{C_\Phi}{n_L(r)}+2^{-n_L(r)}<2^{-2L}r.
\end{equation}
Set
\begin{equation}\label{eq:D}
 \mathcal D_b(r):=\bigcup_{L=1}^{\infty}\ \bigcup_{0\leq K<2^L}
 E_+^{L,K}(n_L(r)).
\end{equation}

\begin{proposition}\label{prop:D}
The set $\mathcal D_b(r)$ is open, $\mu(\mathcal D_b(r))<r$, and every $x\in[0,1]\setminus\mathcal D_b(r)$ satisfies
\begin{equation}\label{eq:D-property}
 \lim_{L\to\infty}\limsup_{N\to\infty}
 \sup_{|a-a'|\leq2^{-L}}
 \frac{\sqrt N D_N^{a,a'}(x)}{\sqrt{2\log\log N}}=0.
\end{equation}
For computable $r$, the set has a uniformly computable enumeration by rational open intervals and computable tail bounds.
\end{proposition}

\begin{proof}
The measure estimate follows from \eqref{eq:nL-choice}:
\[
 \mu(\mathcal D_b(r))
 <\sum_{L=1}^{\infty}2^L2^{-2L}r=r.
\]
An interval of length at most $2^{-L}$ meets at most two adjacent dyadic cells of order $L$. Decomposing its contribution into a bounded number of intervals anchored at dyadic endpoints and applying \eqref{eq:philipp-small-cell} proves \eqref{eq:D-property}. The effective assertion follows because the tail beginning at level $n$ is bounded by $C_\Phi/n+2^{-n}$.
\end{proof}

\begin{proof}[Proof of Theorem~\ref{thm:ulilb}]
The measure assigned to the $\mathcal L$-sets at level $L$ is less than
\[
 b^{2L}\,2^{-(L+1)}b^{-2L}r=2^{-(L+1)}r.
\]
Together with $\mu(\mathcal D_b(r/2))<r/2$, this proves $\mu(U_b^\lil(r))<r$.

Let $x\notin U_b^\lil(r)$. Proposition~\ref{prop:L} gives \eqref{eq:LILa} for every pair of $b$-adic endpoints. Proposition~\ref{prop:D} permits arbitrary intervals to be approximated by $b$-adic intervals: the two endpoint errors have lengths tending to zero and hence vanish after LIL normalization. Consequently,
\[
 \limsup_{N\to\infty}\frac{\sqrt N D_N^b(x)}{\sqrt{2\log\log N}}
 =\sup_{\substack{L\geq1\\0\leq k<k'<b^L}}
 \sigma_{k/b^L,k'/b^L}.
\]
The map $(a,a')\mapsto\sigma_{a,a'}$ is continuous. Indeed, for an arbitrary interval $I=[a,a')$, write $\psi_I=\mathbbm{1}_I-|I|$ and let $P_b$ be the transfer operator of $x\mapsto\{bx\}$. For every mean-zero function $f$ of bounded variation, the one-dimensional Koksma inequality applied to the shifted $b^{-j}$-mesh gives
\[
 \|P_b^j f\|_\infty\leq b^{-j}\operatorname{Var}(f).
\]
Consequently, the covariance series
\[
 \int_0^1\psi_I^2+2\sum_{j\geq1}\int_0^1\psi_I(x)\psi_I(b^jx)\,dx
\]
converges uniformly over all intervals $I$. Each summand is continuous in $(a,a')$ with respect to $L^1$, so the sum is continuous; for $b$-adic endpoints it reduces to \eqref{eq:sigma-aa}. Hence the displayed supremum agrees with $\sup_{0\leq a<a'\leq1}\sigma_{a,a'}$. Fukuyama's explicit computation \cite[Section~5]{Fu08} gives
\[
 \sqrt2\sup_{0\leq a<a'\leq1}\sigma_{a,a'}=L_b,
\]
which proves \eqref{eq:full-LIL-conclusion}.
\end{proof}

\section{Computing Transcendentally LIL-Normal Numbers }\label{sec:alg}

In this final section, we address the computability of transcendentally (LIL)-normal numbers. We first isolate the simple principle behind the algorithms. We then show that finite unions of polynomial preimages are effectively measurable and explain how computable tail bounds allow us to pass to the countable covers constructed in Sections~\ref{sec:const} and~\ref{sec:lil}.

\subsection{The big idea}

The main idea goes back at least to \cite{BF02} in the context of Sierpiński's covering method. Let $U\subset\rr$ be open and suppose that $\mu(U\cap[0,1])\leq1-\varepsilon$ for some $\varepsilon>0$. If the measures of intersections of $U$ with rational intervals can be computed, one may recursively select nested intervals on which a positive measure deficit remains. Their unique limit then lies outside $U$.

\begin{definition}\label{def:meas}
A measurable set $U\subset\rr$ is called \emph{weakly computably measurable on $[0,1]$} (wcm) if there is an algorithm which, given rational numbers $0\leq\alpha<\beta\leq1$ and $n\in\nn$, returns a rational number $l_n(\alpha,\beta)$ satisfying
\begin{equation}\label{eq:wcm}
 \left|\mu(U\cap(\alpha,\beta))-l_n(\alpha,\beta)\right|\leq2^{-n}.
\end{equation}
\end{definition}
Since finite sets have Lebesgue measure zero, the same algorithm also computes the measures of intersections with rational half-open or closed intervals.

\begin{proposition}\label{prop:algo}
Let $U\subset\rr$ be open and wcm on $[0,1]$. Suppose that
\begin{equation}\label{eq:measure-deficit}
 \mu(U\cap[0,1])\leq1-\varepsilon
\end{equation}
for some positive computable number $\varepsilon$. Then there is a computable real number $\nu\in[0,1]\setminus U$. If $\nu$ is irrational, its digits in every integer base $b\geq2$ are computable.
\end{proposition}

\begin{proof}
Replacing $\varepsilon$ by a smaller positive rational number, we may suppose that $\varepsilon\in\mathbb Q$. Put
\[
 I_0=[0,1),\qquad d_0=0,\qquad\varepsilon_0=\varepsilon.
\]
Inductively, suppose that for some $N\geq1$ we have already constructed
\[
 I_{N-1}=\left[\frac{d_{N-1}}{N!},\frac{d_{N-1}+1}{N!}\right)
\]
with
\begin{equation}\label{eq:inductive-deficit}
 \mu(U\cap I_{N-1})\leq\frac1{N!}-\varepsilon_{N-1}.
\end{equation}
Split $I_{N-1}$ into the $N+1$ intervals
\[
 I_{N,k}:=\left[
 \frac{(N+1)d_{N-1}+k}{(N+1)!},
 \frac{(N+1)d_{N-1}+k+1}{(N+1)!}
 \right),\qquad k=0,\ldots,N.
\]
By the wcm property, compute rational numbers $l_{N,k}$ such that
\begin{equation}\label{eq:child-approx}
 \left|l_{N,k}-\mu(U\cap I_{N,k})\right|
 \leq\frac{\varepsilon_{N-1}}{4(N+1)}.
\end{equation}
Choose the smallest $k_N$ for which $l_{N,k}$ is minimal and set
\begin{equation}\label{eq:digit-recursion}
 d_N:=(N+1)d_{N-1}+k_N,
 \qquad I_N:=I_{N,k_N},
 \qquad\varepsilon_N:=\frac{\varepsilon_{N-1}}{2(N+1)}.
\end{equation}
Averaging the approximants and using \eqref{eq:inductive-deficit} gives
\begin{align*}
 l_{N,k_N}
 &\leq\frac1{N+1}\sum_{k=0}^{N}l_{N,k}\\
 &\leq\frac1{(N+1)!}-\frac{\varepsilon_{N-1}}{N+1}
 +\frac{\varepsilon_{N-1}}{4(N+1)}.
\end{align*}
Together with \eqref{eq:child-approx}, this yields
\begin{equation}\label{eq:next-deficit}
 \mu(U\cap I_N)\leq\frac1{(N+1)!}-\varepsilon_N.
\end{equation}
Thus the induction continues.

For each $N$, the interval $I_N$ contains a point $x_N\notin U$ by \eqref{eq:next-deficit}. The intervals are nested and have length $1/(N+1)!$, so their endpoints form a computable Cauchy sequence with a unique limit $\nu$. Since $x_N\to\nu$ and $[0,1]\setminus U$ is closed, one has $\nu\notin U$.

Finally, suppose that $\nu$ is irrational. Given $b,n\in\nn$, compute the nested intervals $I_N$ until one of them is contained in a single interval
\[
 [j b^{-n},(j+1)b^{-n})
\]
for some integer $j$. Such an $N$ exists because $\nu$ is not a $b$-adic rational. The integer $j$ determines the first $n$ base-$b$ digits of $\nu$.
\end{proof}

The remaining task is to verify the wcm property of the covers used above. We begin with finite unions.

\begin{lemma}\label{lem:finite}
Let $V=\bigcup_{k=1}^{K}(a_k,b_k)$ be a finite union of bounded intervals with computable endpoints. Then $V$ is wcm on $[0,1]$.
\end{lemma}

\begin{proof}
If all endpoints are rational, merge the intervals into disjoint components and compute the length of their intersection with any rational interval exactly. In the general case, for a prescribed precision $2^{-n}$ choose rational approximations $\widehat a_k,\widehat b_k$ with $\widehat a_k<\widehat b_k$ and
\[
 |a_k-\widehat a_k|+|b_k-\widehat b_k|\leq K^{-1}2^{-n}.
\]
Then
\[
 \mu\left(
 \left(\bigcup_{k=1}^{K}(a_k,b_k)\right)\mathbin\triangle
 \left(\bigcup_{k=1}^{K}(\widehat a_k,\widehat b_k)\right)
 \right)
 \leq\sum_{k=1}^{K}\bigl(|a_k-\widehat a_k|+|b_k-\widehat b_k|\bigr)
 \leq2^{-n}.
\]
The assertion follows from the rational case and the triangle inequality.
\end{proof}

\subsection{Polynomial preimages}

The next lemma treats polynomial preimages. In particular, it does not assume that the values of a polynomial at the endpoints of a dyadic interval determine its behavior in the interior.

\begin{lemma}\label{lem:polypre}
Let $p_1,\ldots,p_N\in\zz[X]$ and, for each $n$, let $V_n$ be a finite union of bounded open intervals with computable endpoints. Then
\begin{equation}\label{eq:finite-poly-union}
 V:=\bigcup_{n=1}^{N}p_n^{-1}(V_n)
\end{equation}
is wcm on $[0,1]$.
\end{lemma}

\begin{proof}
It is enough to treat a single polynomial $p$ of degree $d$ and a single interval $I=(a,b)$. Fix $m\in\nn$. Choose rational numbers $a<a_m<b_m<b$ such that
\begin{equation}\label{eq:inner-target}
 (a_m-a)+(b-b_m)\leq
 \left(\frac{2^{-m}}{K_d}\right)^d.
\end{equation}
This can be done effectively by successively refining rational enclosures of the computable numbers $a<b$. By Proposition~\ref{prop:poly},
\begin{align}\label{eq:poly-boundary-error}
 &\mu\left([0,1]\cap p^{-1}(I)\setminus
 p^{-1}((a_m,b_m))\right)\\
 &\qquad\leq s_p\bigl((a_m-a)+(b-b_m)\bigr)\leq2^{-m}.
 \nonumber
\end{align}
The set $[0,1]\cap p^{-1}((a_m,b_m))$ is a finite union of intervals whose endpoints belong to $\{0,1\}$ or are real roots of the two integer polynomials obtained by clearing denominators in $p-a_m$ and $p-b_m$. Sturm's algorithm isolates these roots and computes their order. Hence these endpoints are uniformly computable, and Lemma~\ref{lem:finite} applies. Allocating the error among finitely many intervals and polynomials proves the result.
\end{proof}

We now pass to countable unions. The precise formulation keeps separate the tail in the polynomial index and the tail within each target cover.

\begin{proposition}\label{prop:polywcm}
Let $(p_n)_{n\geq1}\subset\zz[X]$ be a computable sequence of non-constant polynomials and let
\begin{equation}\label{eq:countable-poly-cover}
 U:=\bigcup_{n=1}^{\infty}p_n^{-1}(U_n),
 \qquad U_n:=\bigcup_{k=1}^{\infty}I_{n,k},
\end{equation}
where the $I_{n,k}$ are bounded open intervals with uniformly computable endpoints. Assume the following.
\begin{enumerate}
\item There is a computable sequence $s_m$ such that
\begin{equation}\label{eq:outer-tail}
 \mu\left([0,1]\cap\bigcup_{n=s_m}^{\infty}p_n^{-1}(U_n)\right)
 \leq2^{-m}.
\end{equation}
\item For every $n$ there is a computable sequence $t_{n,m}$ such that
\begin{equation}\label{eq:inner-tail}
 \mu\left(\bigcup_{k=t_{n,m}}^{\infty}I_{n,k}\right)\leq2^{-m}.
\end{equation}
\end{enumerate}
Then $U$ is wcm on $[0,1]$.
\end{proposition}

\begin{proof}
Fix a desired precision $2^{-m}$ and put $s:=s_{m+2}$. The contribution of the indices $n\geq s$ is at most $2^{-m-2}$. For each $1\leq n<s$, choose an integer $j_n$ so large that
\begin{equation}\label{eq:choose-jn}
 K_{d_{p_n}}2^{-j_n/d_{p_n}}
 \leq\frac{2^{-m-2}}{\max\{1,s-1\}},
\end{equation}
and set $t_n:=t_{n,j_n}$. Proposition~\ref{prop:poly} and \eqref{eq:inner-tail} imply
\begin{align*}
 &\mu\left([0,1]\cap\bigcup_{n=1}^{s-1}
 p_n^{-1}\left(\bigcup_{k=t_n}^{\infty}I_{n,k}\right)\right)\\
 &\qquad\leq\sum_{n=1}^{s-1}K_{d_{p_n}}2^{-j_n/d_{p_n}}
 \leq2^{-m-2}.
\end{align*}
The remaining set
\[
 U^{(m)}:=\bigcup_{n=1}^{s-1}p_n^{-1}
 \left(\bigcup_{k=1}^{t_n-1}I_{n,k}\right)
\]
is wcm by Lemma~\ref{lem:polypre}. Compute its intersection measure with error at most $2^{-m-1}$. The two tail estimates and the triangle inequality give a total error at most $2^{-m}$.
\end{proof}

\subsection{Putting the pieces together}

We first treat ordinary transcendental normality. Recall the Sierpiński cover $U(r)$ from \eqref{eq:sircover}. The following quantitative truncation is proved in Appendix~\ref{sec:app1}.

\begin{lemma}\label{lem:trunc}
For every rational $r\in(0,1)$ and every $k\in\nn$ one can compute a finite union $U^{[k]}(r)$ of intervals with rational endpoints such that
\begin{equation}\label{eq:sierpinski-truncation}
 U^{[k]}(r)\subset U(r),\qquad
 \mu\bigl(U(r)\setminus U^{[k]}(r)\bigr)\leq2^{-k}.
\end{equation}
For instance, one may retain only
\begin{equation}\label{eq:explicit-truncation-indices}
\begin{gathered}
 2\leq b\leq B_k:=2^{k+4},\qquad
 1\leq m\leq M_k:=2^{k+4},\\
 n_{m,b}(r)\leq n\leq N_k:=12\cdot2^{7k+32}.
\end{gathered}
\end{equation}
\end{lemma}

The same conclusion holds for $U^+(r)$. Indeed, first truncate the integer translations in
\[
 U^+(r)=\bigcup_{z\in\zz}\left(z+U\left(\frac r{2^{|z|+2}}\right)\right);
\]
the omitted translates have total measure bounded by $2^{-|z_0|-1}r$ when $|z|>|z_0|$. Then apply Lemma~\ref{lem:trunc} to the finitely many remaining covers.

Fix a rational $r\in(0,1)$ and write
\[
 \delta_n:=\left(\frac{2^{-n}r}{K_{d_{p_n}}}\right)^{d_{p_n}}.
\]
The polynomial tail of $U_{\mathcal P}(r)$ satisfies
\begin{equation}\label{eq:normal-outer-tail}
 \mu\left([0,1]\cap\bigcup_{n=s}^{\infty}
 p_n^{-1}(U^+(\delta_n))\right)
 <r\sum_{n=s}^{\infty}2^{-n}\leq2^{1-s}.
\end{equation}
Enumerating the intervals in successive expanding boxes as in Lemma~\ref{lem:trunc}, and listing every interval only once, turns the same estimates into the tail bound required in the second hypothesis of Proposition~\ref{prop:polywcm}. Hence $U_{\mathcal P}(r)$ is wcm on $[0,1]$. Taking, for example, $r=1/2$, Proposition~\ref{prop:algo} produces a computable number
\[
 \nu\in[0,1]\setminus U_{\mathcal P}(1/2).
\]
It is transcendentally normal and therefore irrational, so its digits are computable in every base. This completes the proof of Theorem~\ref{thm:main}.

We finally consider LIL-normality.

\begin{lemma}\label{lem:liltrunc}
Let $r\in(0,1)$ be computable. There is a uniformly computable enumeration
\begin{equation}\label{eq:lil-interval-enumeration}
 U^\lil(r)=\bigcup_{j=1}^{\infty}I_j(r)
\end{equation}
by bounded open intervals with rational endpoints and a computable sequence $t_k(r)$ such that
\begin{equation}\label{eq:lil-finite-approximation}
 \mu\left(\bigcup_{j=t_k(r)}^{\infty}I_j(r)\right)\leq2^{-k}.
\end{equation}
The construction is uniform in $r$.
\end{lemma}

\begin{proof}
We spell out the truncation because it is needed for the algorithm. Enumerate the constituent exceptional events in successive finite stages, always listing all interval components of an event consecutively and never listing the same event twice. In the base decomposition \eqref{eq:all-bases-cover}, the tail over $b>B$ has measure at most
\[
 r\sum_{b>B}2^{-b}=r2^{-B}.
\]
Thus only finitely many bases are needed at any prescribed accuracy.

Fix one of these bases and consider $U_b^\lil(s)$ from \eqref{eq:covUb}. For the set $\mathcal D_b(s/2)$, the total measure assigned to dyadic levels larger than $L_0$ is bounded by a geometric tail. At a fixed level $(L,K)$, truncating the union defining $E_+^{L,K}(n_L(s/2))$ at a level $n_0$ leaves a set of measure at most
\[
 \frac{C_\Phi}{n_0}+2^{-n_0}
\]
by Lemma~\ref{lem:phi}. Both bounds tend effectively to zero.

For the second part of \eqref{eq:covUb}, the total measure assigned to orders $L>L_0$ is again bounded by a geometric tail. For each of the finitely many retained pairs of $b$-adic endpoints, the tail over $M>M_0$ in \eqref{eq:L} is at most a geometric sum of the budgets $2^{-M}$. For fixed $M$, truncating the unions in \eqref{eq:Yup}--\eqref{eq:Zup} at $q=q_0$ leaves measure bounded by the corresponding computable functions in \eqref{eq:Yup-bound}--\eqref{eq:Zup-bound}, with $N$ replaced by $q_0+1$. These bounds tend effectively to zero by Proposition~\ref{prop:RW}.

Allocate the error $2^{-k}$ among these tails and choose every cutoff at stage $k+1$ at least as large as at stage $k$. Every retained event is a finite union of base-$b$ cylinder intervals enlarged by a rational radius, and hence is a finite union of rational open intervals. Let $t_k(r)$ be the first index after all events retained at stage $k$ have been listed. By the preceding union bounds, every event listed from that index onward belongs to one of the omitted tails, so \eqref{eq:lil-finite-approximation} holds. The union of all stages contains every defining event and therefore equals $U^\lil(r)$. All choices are effective because the tail bounds and cutoffs in Section~\ref{sec:lil} are computable.
\end{proof}

As before, truncating the integer translations and applying Lemma~\ref{lem:liltrunc} gives a computable interval enumeration of $U_+^\lil(r)$ with computable measure tails. The outer polynomial tail of $U_{\mathcal P}^\lil(r)$ is bounded exactly as in \eqref{eq:normal-outer-tail}. Proposition~\ref{prop:polywcm} therefore shows that $U_{\mathcal P}^\lil(r)$ is wcm. Taking $r=1/2$ in Proposition~\ref{prop:algo} produces a computable
\[
 \nu_\lil\in[0,1]\setminus U_{\mathcal P}^\lil(1/2).
\]
This number is transcendentally LIL-normal, hence irrational, and its digits can be computed in every integer base. Theorem~\ref{thm:main2} follows.

\appendix
\section{Quantitative Bounds for Sierpiński's Cover}\label{sec:app1}

We give the quantitative estimates used in Lemma~\ref{lem:trunc}. Recall the sets $U_{b,m,n,d}$ from \eqref{eq:sirint}.

\begin{lemma}\label{lem:intbound}
For integers $b\geq2$, $m,n\geq1$, and $0\leq d<b$,
\begin{equation}\label{eq:intbound}
 \mu(U_{b,m,n,d})\leq\frac{12m^4}{bn^2}.
\end{equation}
\end{lemma}

\begin{proof}
Each interval in \eqref{eq:sirint} has length $3b^{-n}$. Interpreting a word $(q_1,\ldots,q_n)$ as $n$ independent uniform base-$b$ digits, the counting variable $C_d$ is binomial with success probability $1/b$. Hence
\begin{align*}
 \mu(U_{b,m,n,d})
 &\leq3\,\pp\left(\left|\frac{C_d}{n}-\frac1b\right|\geq\frac1m\right).
\end{align*}
Let $X$ be a Bernoulli variable with parameter $1/b$, put $\sigma^2=\operatorname{Var}(X)$ and $\mu_4=\ee(X-1/b)^4$. Independence gives the exact fourth-moment identity
\begin{equation}\label{eq:binomial-fourth}
 \ee\left(C_d-\frac nb\right)^4
 =n\mu_4+3n(n-1)\sigma^4.
\end{equation}
Since $\mu_4\leq1/b$ and $\sigma^4\leq1/b^2\leq1/b$, it follows that
\[
 \ee\left(\frac{C_d}{n}-\frac1b\right)^4\leq\frac4{bn^2}.
\]
Markov's inequality proves \eqref{eq:intbound}.
\end{proof}

\begin{lemma}\label{lem:seriesbounds}
For $b,N\geq2$ and $m\geq1$,
\begin{equation}\label{eq:tail-n}
 \mu\left(\bigcup_{n=N}^{\infty}\bigcup_{d=0}^{b-1}U_{b,m,n,d}\right)
 \leq\frac{12m^4}{N-1}.
\end{equation}
Moreover, for $b\geq2$ and $M\geq2$,
\begin{equation}\label{eq:tail-m}
 \mu\left(\bigcup_{m=M}^{\infty}\bigcup_{n=n_{m,b}(r)}^{\infty}
 \bigcup_{d=0}^{b-1}U_{b,m,n,d}\right)
 \leq\frac{r}{2b^2(M-1)},
\end{equation}
and for $B\geq2$,
\begin{equation}\label{eq:tail-b}
 \mu\left(\bigcup_{b=B}^{\infty}\bigcup_{m=1}^{\infty}
 \bigcup_{n=n_{m,b}(r)}^{\infty}\bigcup_{d=0}^{b-1}U_{b,m,n,d}\right)
 <\frac r{B-1}.
\end{equation}
\end{lemma}

\begin{proof}
By Lemma~\ref{lem:intbound} and the telescoping estimate
\[
 \sum_{n=N}^{\infty}\frac1{n^2}
 \leq\sum_{n=N}^{\infty}\frac1{n(n-1)}=\frac1{N-1},
\]
we obtain \eqref{eq:tail-n}. Since
\[
 n_{m,b}(r)-1\geq\frac{24m^6b^2}{r},
\]
formula \eqref{eq:tail-n} gives
\[
 \mu\left(\bigcup_{n=n_{m,b}(r)}^{\infty}\bigcup_{d=0}^{b-1}U_{b,m,n,d}\right)
 \leq\frac r{2m^2b^2}.
\]
Summing over $m\geq M$ proves \eqref{eq:tail-m}; summing first over all $m$ and then over $b\geq B$ proves \eqref{eq:tail-b}.
\end{proof}

\begin{proof}[Proof of Lemma~\ref{lem:trunc}]
Let $B_k,M_k,N_k$ be as in \eqref{eq:explicit-truncation-indices} and retain exactly the corresponding finite collection of intervals. The omitted part is contained in the union of three tails: $b>B_k$, or $b\leq B_k$ and $m>M_k$, or $b\leq B_k$, $m\leq M_k$ and $n>N_k$. Lemma~\ref{lem:seriesbounds} gives
\begin{align*}
 T_b&<\frac r{B_k},\\
 T_m&\leq\sum_{b=2}^{B_k}\frac r{2b^2M_k}
 \leq\frac1{2M_k},\\
 T_n&\leq\frac{12B_k}{N_k}
 \sum_{m=1}^{M_k}m^4
 \leq\frac{12B_kM_k^5}{N_k}.
\end{align*}
With the explicit choices in \eqref{eq:explicit-truncation-indices},
\[
 T_b+T_m+T_n<2^{-k}.
\]
All retained endpoints are rational and are computable uniformly from the rational parameter $r$.
\end{proof}

\section{Elementary proof of a weak P\'{o}lya-type inequality}\label{sec:polya}

For completeness, we give an elementary bound which is weaker than Lemma~\ref{lem:polya}, but still sufficient for all computability arguments in this paper.

\begin{lemma}\label{lem:matrix}
For $d\geq1$, let $A^{(d)}$ be the $(d+1)\times(d+1)$ Vandermonde matrix
\[
 A^{(d)}_{i,j}=i^{j-1},\qquad1\leq i,j\leq d+1.
\]
Its least singular value satisfies
\begin{equation}\label{eq:singular-lower}
 s_{\min}(A^{(d)})\geq(d+1)^{-d(d+1)}=:C_d^{-1}.
\end{equation}
\end{lemma}

\begin{proof}
The Vandermonde determinant is nonzero, so $A^{(d)}$ is invertible. More explicitly, if
\[
 \sum_{j=0}^{d}\gamma_{j+1}i^j=0,
 \qquad i=1,\ldots,d+1,
\]
then the polynomial $\sum_{j=0}^{d}\gamma_{j+1}x^j$ has $d+1$ distinct zeros and hence all coefficients vanish.

The product of the singular values equals $|\det A^{(d)}|$, which is a positive integer and therefore at least one. Every singular value is at most the Frobenius norm, and
\[
 \|A^{(d)}\|_F^2
 \leq(d+1)^2(d+1)^{2d}=(d+1)^{2d+2}.
\]
Consequently,
\[
 s_{\min}(A^{(d)})
 \geq\|A^{(d)}\|_F^{-d}
 \geq(d+1)^{-d(d+1)}.
\]
\end{proof}

\begin{lemma}\label{lem:small}
Let $q\in\zz[X]+\rr$ be a polynomial of degree $d>0$ with integer coefficients except for a possible real constant term. Then, for every $\varepsilon>0$,
\begin{equation}\label{eq:weak-polya}
 \mu\{x\in[0,1]:|q(x)|\leq\varepsilon^d\}
 \leq4d^3(d+1)^{d+1}\varepsilon.
\end{equation}
\end{lemma}

\begin{proof}
The case $d=1$ is immediate, so assume $d\geq2$. Put
\[
 E:=\{x\in[0,1]:|q(x)|\leq\varepsilon^d\}.
\]
The boundary points of $E$ in $(0,1)$ are zeros of the polynomial $q^2-\varepsilon^{2d}$, which has degree $2d$. Hence $E$ has at most $d+1$ connected components.

Fix $y\in\rr$ and expand
\[
 q(y+t)=\sum_{j=0}^{d}b_jt^j.
\]
For $\delta>0$, set
\[
 w_i:=q(y+i\delta),\qquad
 v:=(b_0,b_1\delta,\ldots,b_d\delta^d)^T,
 \qquad i=1,\ldots,d+1.
\]
Then $w=A^{(d)}v$. Since the leading coefficient $b_d$ is a nonzero integer, Lemma~\ref{lem:matrix} gives
\begin{equation}\label{eq:grid-value}
 \max_{1\leq i\leq d+1}|q(y+i\delta)|
 \geq\frac{\delta^d}{C_d\sqrt{d+1}}
 =\left(\frac\delta{L_d}\right)^d,
 \qquad L_d:=(C_d\sqrt{d+1})^{1/d}.
\end{equation}
Choose $\delta=2L_d\varepsilon$. If a connected component of $E$ had length greater than $(d+1)\delta$, one could choose $y$ so that all points $y+\delta,\ldots,y+(d+1)\delta$ lie in that component. This would contradict \eqref{eq:grid-value}. Therefore
\[
 \mu(E)\leq2(d+1)^2L_d\varepsilon.
\]
Since $L_d=(d+1)^{d+1+1/(2d)}$ and $d\geq2$, the elementary inequality
\[
 2(d+1)^2L_d\leq4d^3(d+1)^{d+1}
\]
proves \eqref{eq:weak-polya}.
\end{proof}

\section*{Acknowledgments}

The author would like thank  Christoph Aistleitner for many helpful comments on an earlier draft and Verónica Becher for spotting a mistake in the proof.  Special thanks go to the anonymous referee for carefully checking this manuscript and pointing out several flaws in the proof. This work was funded by the Deutsche Forschungsgemeinschaft (DFG, German Research Foundation) - 558731723.

\section*{AI Declaration}
The original draft has been crafted without any use of AI. In the revised version, an LLM has been used to correct and streamline some parts of Section~\ref{sec:lil}: the Berry-Essen argument in the proof of Proposition~\ref{prop:RW}, a simplified definition of the covers in Section~\ref{sec:bcover} and a more explicit version of Lemma~\ref{lem:philipp-decomp}. Apart from that, AI has been only used for minor language and typographical corrections.


\begin{thebibliography}{55}
 \bibitem{AB13} C. Aistleitner, I. Berkes. \newblock Limit distributions in metric discrepancy
theory. \newblock {\em  Monatsh. Math.}, 169(3-4):253–265 (2013).
\bibitem{ABSS17} C.Aistleitner, V. Becher, A.-M. Scheerer, T. Slaman. \newblock On the construction of absolutely normal numbers. \newblock {\em Acta Arithmetica} 180, 333-346 (2017).

\bibitem{AFP24} C. Aistleitner, L. Frühwirth, J. Prochno. \newblock  Diophantine conditions in the law of the iterated logarithm for lacunary systems. \newblock {\em Probab. Theory Relat. Fields} 192, 545–574 (2025).

\bibitem{ABS22} A. Algom, S. Baker, P. Shmerkin. \newblock
On normal numbers and self-similar measures. \newblock 
{\em Adv. Math.} 399, 108276 (2022).

\bibitem{AB11} N. Alvarez, V. Becher. \newblock M. Levin's construction of absolutely normal numbers with very low discrepancy. \newblock {\em Mathematics of Computation } 86 (308), 2927 – 2946 (2017).

\bibitem{AL23} A. Aveni, P. Leonetti. \newblock Most numbers are not normal. \newblock { \em Mathematical Proceedings of the Cambridge Philosophical Society}, 175(1), pp. 1–11 (2023).

\bibitem{BBCP04} D.H. Bailey, J.M. Borwein, R.E. Crandall, C. Pomerance. On the binary expansions of algebraic numbers. \newblock {\em J. Th\'{e}or. Nombres Bordeaux} 16(3), 487-518 (2004).

\bibitem{BB25}
S. Baker and A. Banaji. Polynomial Fourier decay for fractal measures and their pushforwards. \emph{Math. Ann.} 392, no. 1, 209--261 (2025).

\bibitem{BBS16} V. Becher, Y. Bugeaud, T.A. Slaman. \newblock On Simply Normal Numbers to Different Bases. {\em 
Mathematische Annalen} 364 (1), 125-150 (2016).
\bibitem{BF02}  V. Becher, S. Figueira. \newblock An example of a computable absolutely normal number. \newblock {\em Theoretical Computer Science} 270,  947–958 (2002).
\bibitem{BFP07} V. Becher, S.Figueira, R. Picchi. \newblock Turing’s unpublished algorithm for normal numbers. \newblock {\em Theoretical Computer Science}, 377:126–138 (2007).
 \bibitem{BHS13} V. Becher, P. A. Heiber, T.A. Slaman. \newblock A polynomial-time algorithm for computing absolutely normal
numbers. \newblock {\em Information and Computation} 232, 1-9 (2013).

\bibitem{BM22} V. Becher, M. G. Madritsch. \newblock On a question of Mendes France on normal numbers. \newblock {\em Acta Arithmetica} 203, 271-288 (2022).
\bibitem{BY19} V. Becher, S. Yuhjtman. \newblock On absolutely normal and continued fraction normal numbers. \newblock {\em International
Mathematics Research Notices} 19 (2019).
\bibitem{Ber76} I. Berkes. \newblock On the asymptotic behaviour of $\sum f(n_kx)$. I. Main Theorems, II. Applications. \newblock {\em Z. Wahr. verw. Geb.} 34, 319-345, 347-365 (1976).
\bibitem{Birk31} G. D. Birkhoff. \newblock Proof of the Ergodic Theorem. {\em Proc. Nat. Acad. Sci.} 17, 656–660 (1931).
\bibitem{Bor09} E. Borel. \newblock  Les probabilit\'{e}s d\'{e}nombrables et leurs applications arithm\'{e}tiques.\newblock {\em Supplemento di Rendiconti
del Circolo Matematico di Palermo} 27 , 247–271 (1909).
\bibitem{Bug12} Y. Bugeaud.\newblock On the expansions of a real number to several integer bases. \newblock {\em Revista
Matemática Iberoamericana}, 28(4):931–946 (2012).

\bibitem{Cas59} J.W.S Cassels. \newblock On a problem of Steinhaus about normal numbers. {\em Colloq. Math. 7}, 95–101 (1959).

\bibitem{DGW24}
Y. Dayan, A. Ganguly and B. Weiss. \newblock
Random walks on tori and normal numbers in self-similar sets.
{\em Amer. J. Math.} 146, 467--493 (2024).

\bibitem{Egg49} H. G. Egglestone. \newblock  The fractional dimension of a set defined by decimal properties. \newblock {\em Quart. J. Math.},
Oxford Ser. 20, 31–36 (1949).
\bibitem{Fu08} K. Fukuyama. \newblock The law of the iterated logarithm for discrepancies of $\{\theta^nx\}$. \newblock { \em Acta Mathematica Hungarica}, 118(1):155–170 (2008).
\bibitem{GG64} S. G\'{a}l, L. G\'{a}l. \newblock The discrepancy of the sequence $\{(2^nx)\}$. \newblock {\em Koninklijke
Nederlandse Akademie van Wetenschappen Proceedings. Seres A 67 = Indagationes
Mathematicae} , 26:129–143 (1964).

\bibitem{HS15}
M. Hochman, P. Shmerkin.
Equidistribution from fractal measures.
{\em Invent. Math.} 202, 427--479 (2015).

\bibitem{Knu97} D.E. Knuth. The Art of Computer Programming. Volume 2: Seminumerical Algorithms. Third Edition. Reading, Massachusetts: Addison-Wesley, 1997.

\bibitem{Leb17} H. Lebesgue. \newblock  Sur certaines d\'{e}monstrations d’existence. {\em Bulletin de la Soci\'{e}t\'{e}
Math\'{e}matique de France} 45:132–144 (1917).


\bibitem{Lev79} M. Levin. \newblock On absolutely normal numbers. {\em Vestnik Moskovskogo Universiteta.
Seriya 1. Matematika. Mekhanika}, 1:31–37, 87, 1979. English translation in Moscow University Mathematics Bulletin, 34  no. 1, 32-39 (1979).

\bibitem{LM21} J. H. Lutz, E. Mayordomo. \newblock Computing absolutely normal numbers in nearly linear time. \newblock {\em Information and Computation} 281, 104746 (2021).

\bibitem{Mad18} M. Madritsch. \newblock Normal numbers and symbolic dynamics. \newblock {\em Sequences, groups, and number theory}, pp. 271–329 (2018).

\bibitem{Man25} C. Manai. \newblock 
Super-dense sets and their role in the theory of normal numbers
\emph{Adv. Math.} 503, 111205 (2026).

\bibitem{Man26} C. Manai. \newblock Digit Mixing under Polynomial Maps \newblock {\em Preprint } 	arXiv:2606.08325 

\bibitem{Man26b} C. Manai. \newblock On Normality Preserving Operations {\em Preprint } 	arXiv:2609.24665 

\bibitem{Ol04a} L. Olsen. \newblock Applications of multifractal divergence points to sets of numbers defined by their N-adic expansion. \newblock {\em  Math. Proc. Camb. Phil. Soc. 136}, no. 1, 139–165 (2004).

\bibitem{Ol04b} L. Olsen. \newblock Extremely non-normal numbers. \newblock {\em Math. Proc. Camb. Phil. Soc. 137}, no. 1, 43–53 (2004).

\bibitem{Phi75} W. Philipp. \newblock Limit theorems for lacunary series and uniform distribution mod 1. \newblock {\em Acta Arithmetica}, 26(3):241–251 (1975).

\bibitem{Pol28} G. P\'{o}lya. Beitrag zur Verallgemeinerung des Verzerrungssatzes auf mehrfach zusammenhängende Gebiete. { \em Sitzungsberichte Akad. Berlin:} 280–282 (1928).

\bibitem{Rau76} G. Rauzy. \newblock Nombres normaux et processus d\'{e}terministes. \newblock {\em
Acta Arithmetica}, 29:211–225 (1976).

\bibitem{Rem36} E.J. Remez. Sur une propri\'{e}t\'{e} des polyn\^{o}mes de Tchebyscheff. {\em Comm. Inst. Sci. Kharkow.} 13: 93–95 (1936).

\bibitem{Riv08} T. Rivoal. \newblock On the bits counting function of real numbers. \newblock {\em Journal Australian Mathematical Society} 85, 95–111 (2008).

\bibitem{Schee17} A.-M. Scheerer. \newblock On the continued fraction expansion of absolutely normal numbers.
\newblock {\em Preprint, arXiv:1701.07979v1} (2017).

\bibitem{Schm62} W. M. Schmidt. \newblock Über die Normalität von Zahlen zu verschiedenen Basen. \newblock  {\em Acta Arithmetica}, 7:299–309 (1961/1962).

 \bibitem{Sir17} M.W. Sierpiński. \newblock  D\'{e}monstration \'{e}l\'{e}mentaire du th\'{e}or\`{e}me de M. Borel sur les nombres absolument
normaux et d\'{e}termination effective d’un tel nombre. \newblock {\em Bull. Soc. Math. France} 45, 127–132 (1917).

\bibitem{Str64} V. Strassen. \newblock An invariance principle for the law of iterated logarithm. \newblock {\em Z. Wahr. verw. Geb.} 3, 211-226 (1964).

\bibitem{Str67} V. Strassen. \newblock Almost sure behavior of sums of independent random variables and martingales. {\em Proc. 5-th Berkeley Sympos. Math. Statist. Probab.,} Vol. II (Part I) 315-343. Univ. of Calif. Press, Berkeley (1967).

\bibitem{Tur92}  A. M. Turing. \newblock A note on normal numbers. \newblock { \em in: J.L. Britton (Ed.), Collected Works of A.M. Turing: Pure Mathematics. } North Holland,
Amsterdam, 1992, pp. 117–119.

\bibitem{Van18} J. Vandehey. \newblock On the binary digits of $\sqrt{2}$. {\em Integers } 18 (2018).

\bibitem{YD97} Y. Yi-Sheng, G. Dun-he. \newblock A note on a lower bound for the smallest singular value. Linear Algebra and its Applications 253 (1-3), pp. 25-38 (1997).

\bibitem{Kac46} M. Kac. \newblock On the distribution of values of sums of the type $\sum f(2^kt)$. \newblock {\em Ann. of Math. (2)} 47, 33--49 (1946).

\bibitem{Mar50} G. Maruyama. \newblock On an asymptotic property of a gap sequence. \newblock {\em K\^{o}dai Math. Sem. Rep.} 2, 31--32 (1950).

\bibitem{Pet75} V. V. Petrov. \newblock {\em Sums of Independent Random Variables}. \newblock Springer, Berlin (1975).

\end{thebibliography}
\end{document}